\documentclass[conf]{new-aiaa}
\usepackage[utf8]{inputenc}

\usepackage{bm}
\usepackage{amsmath}
\usepackage{subfig}
\usepackage{xcolor}

\usepackage{graphicx}
\usepackage[version=4]{mhchem}
\usepackage{siunitx}
\usepackage{longtable,tabularx}
\usepackage[preprint]{}

\usepackage{enumitem}
\usepackage[mathscr]{eucal}

\DeclareMathOperator{\oec}{\mathbf{o\mkern-3mu e}}

\DeclareMathOperator{\doe}{\delta\mkern-3mu\oec}

\title{Spacecraft Relative Motion Dynamics and Control Using Fundamental Solution Constants}

\author{Ethan R. Burnett\footnote{Ph.D., Ann and H.J. Smead Department of Aerospace Engineering Sciences, University of Colorado, 431 UCB, Colorado Center for Astrodynamics Research, Boulder, CO 80309} and Hanspeter Schaub\footnote{Professor, Glenn L. Murphy Chair of Engineering, Smead Department of Aerospace Engineering Sciences, University of Colorado, 431 UCB, Colorado Center for Astrodynamics Research, Boulder, CO 80309-0431. AAS Fellow, AIAA Fellow.}}
\affil{Ann and H.J. Smead Department of Aerospace Engineering Sciences, University of Colorado Boulder}

\begin{document}

\maketitle

\begin{abstract}
This paper explores expressing the relative state in the close-proximity satellite relative motion problem in terms of fundamental solution constants. The nominal uncontrolled relative state can be expressed in terms of a weighted sum of fundamental and geometrically insightful motions. These fundamental motions are obtained using Lyapunov-Floquet theory. In the case that the dynamics are perturbed by the action of a controller or by unmodeled dynamics, the weights on each fundamental solution are allowed to vary as in a variation-of-parameters approach, and in this manner function as state variables. This methodology reveals interesting insights about satellite relative motion and also enables elegant control approaches. This approach can be applied in any dynamical environment as long as the chief orbit is periodic, and this is demonstrated with results for relative motion analysis and control in the eccentric Keplerian problem and in the circular restricted three-body problem (CR3BP). Some commentary on extension of the methodology beyond the periodic chief orbit case is also provided. This is a promising and widely applicable new approach to the close-proximity satellite relative motion problem. \footnotetext{\\ Copyright \textcopyright \ 2021 Ethan R. Burnett and Hanspeter Schaub}
\end{abstract}

\section{Introduction}
In the close-proximity relative motion control problem, the choice of coordinate representation greatly influences the ease of relative motion planning and design, and can also affect the complexity of feedback control development \cite{schaub}. The traditional local relative Cartesian and spherical coordinate representations have proven to be quite popular \cite{ButcherBurnett2017},  but there are other formulations that enable geometric approaches or provide great geometric insight into the influence of perturbations and control. Orbit element differences \cite{RelativeOrbitPaper} and the related relative orbit element (ROE) formulations are important examples \cite{DAmico_EccInc,Guffanti_ICs}, and some missions have already implemented ROE formulations for relative motion control \cite{Amico:2011ff}. As opposed to a local coordinate formulation, depending on the choice of quantities, most or all orbit element differences or ROEs will be constants for uncontrolled Keplerian relative motion. In this sense, the unperturbed dynamics of relative motion are trivial in the orbit element differences, and the mapping from ROEs to other representations of the relative state can provide a local geometric interpretation \cite{LocalCart}. These formulations also generally linearize quite well \cite{SinclairJGCD}, which is a very useful feature.

Another relevant representation of relative motion is discussed in Reference~\citenum{Bennett2016c}, which focuses on the Linearized Relative Orbit Elements (LROEs) for the Clohessy-Wiltshire (CW) problem \cite{clohessy1}, and seeks a representation with a straightforward geometric interpretation. This is accomplished by representing the CW solution in terms of fundamental amplitude, phase, and offset quantities, which are themselves nonlinear functions of the initial relative state quantities in local Cartesian coordinates. These quantities are thus constant in the unperturbed problem in the absence of control. Given nonzero perturbations or control, these quantities vary in time, with their behavior given by equations obtained using a Lagrange-Bracket variational methodology \cite{schaub}. There has also been some work relating integration constants for the CW problem to the ROE formulation. For example, Reference~\citenum{Liou20111865} provides a mapping between the CW integration constants and the general ROEs.

Because the LROEs discussed in Reference~\citenum{Bennett2016c} are developed for CW dynamics, they are only appropriate for near-circular orbits, and would not be appropriate for eccentric orbits. Additionally, the orbit element difference and ROE sets suffer from their own shortcomings. 
For example, the local out-of-plane dynamic decoupling observed in the Keplerian satellite relative motion problem is commonly obscured by formulations using orbit elements \cite{LocalCart}. They are also not directly related to system observables \cite{SullivanDAmico_Survey}. Additionally, the application of any of these integral representations in highly non-Keplerian settings (such as in the three-body problem) would be cumbersome, because while they can still be a valid state representation via variation-of-parameters, they are no longer nearly constant. 

\begin{figure}[h!]
\centering
\includegraphics[]{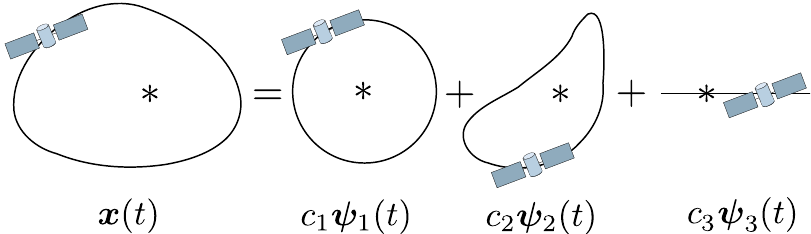}
\caption{Example Satellite Relative Motion as a Sum of Individual Modes} \vspace{-5mm}
\label{fig:ModalConc}
\end{figure}

Motivated by the many benefits of integral relative state representations such as the ROEs and LROEs, and also by their shortcomings, this paper discusses an approach to modeling satellite relative motion dynamics and control that can be applied to the CW problem, the Keplerian problem with eccentricities $0 < e < 1$, and beyond, to the restricted three-body problems. This approach is called the \textit{method of fundamental modal solutions}. The main idea of this work is that satellite relative motion can be represented as a weighted sum of fundamental solutions, chosen for maximal ease of geometric interpretation. 
Figure \ref{fig:ModalConc} illustrates this conceptually with a depiction of bounded planar relative motion decomposed into three simpler constituent modal motions. 
This is accomplished by using the modal solutions previously developed by the authors through an application of Lyapunov-Floquet theory \cite{Burnett_AA2021}, as well as more recent extensions. A benefit of this approach is that stable, unstable, oscillatory, and drift motions are naturally separated. In the absence of perturbation or control effort, the weights on each fundamental solution are constant. However, in the perturbed or controlled problem, the constants are made to vary through a variation-of-parameters framework, such that the weighted sum of the fundamental solutions still always instantaneously describes the relative state. The collection of six time-varying weights serves as the state vector for the problem, and the benefits of this approach for dynamics, visualization, and control are discussed in this paper. This procedure is presently developed for close-proximity relative motion in the vicinity of any closed (periodic) chief orbit, regardless of the governing dynamics. To demonstrate its breadth of potential uses, the methodology is applied to the relative motion problem for Keplerian orbits of any eccentricity, and to orbits in the circular restricted three-body problem (CR3BP). 

There are a few works that are particularly relevant to this paper due to their shared methods, mathematical techniques, and perspectives.  First, this paper employs Lyapunov-Floquet theory \cite{Nayfeh:1979uq}, whose application to the relative motion problem has seen limited but noteworthy study. In Reference~\citenum{Sherrill2015LyapunovFloquet}, a Lyapunov-Floquet (LF) transformation relating the Clohessy-Wiltshire and Tschauner Hempel dynamics is exploited for relative motion control design. 
Additionally, a simple LF transformation is used in Reference~\citenum{KoenigJ2} when incorporating the secular effect of the $J_{2}$ perturbation on orbit element differences. Recently, Reference~\citenum{Gaias_ROEs2021} devises a unified ROE-CW framework for close-proximity satellite relative motion control, adapting work from Reference~\citenum{Liou20111865} and a periodic transformation between ROEs and local Cartesian coordinates for the CW problem. Also, the Keplerian modal decomposition discussed in this paper was first derived analytically in a recently published paper by the author \cite{Burnett_AA2021}. Lastly, other works have previously explored formation flying in the Earth-Moon three-body problem, particularly in the vicinity of the libration points. Past works such as References~\citenum{GomezHowell} and \citenum{HowellL1L2} and references therein are relevant.

This paper is organized as follows. First, the concept of the relative motion modal decomposition is introduced, and there is also a brief review of the Lyapunov-Floquet transformation, which is used extensively in this work. The modal decomposition for general Keplerian orbits is introduced, applicable to any eccentricity $0<e<1$. Then, the method of variation-of-parameters is used to derive the dynamics of the modal constants for the case that the chief orbit is perturbed from its reference periodic motion. Afterwards, the theory is applied to relative motion in the vicinity of periodic orbits in the CR3BP, and an interesting parallel between the Keplerian close-proximity relative motion problem and its CR3BP counterpart is highlighted. Then, applications of the fundamental solution constants are discussed for modeling perturbed relative motion and for efficiently computing relative motion control policies. Finally, numerical results are shown for both the case of Keplerian and CR3BP dynamics. For the Keplerian case, the relative motion modes are computed and used to explore the design space of relative orbits, and the effect of $J_{2}$ on the Keplerian modes is explored. Then, fuel-optimal impulsive relative motion control using the modal decomposition is explored for both Keplerian relative motion and relative motion in the CR3BP. The paper ends with a discussion of future work that would generalize this method to efficient relative motion analysis, modeling, and control in a wide variety of environments.
\section{The Relative Motion Modal Decomposition}
\subsection{Theory and Fundamentals}
In this work, close-proximity satellite relative motion is represented as a weighted sum of fundamental solutions:
\begin{equation}
\label{FS_rm1}
\bm{x}(t) = \sum_{i = 1}^{6}c_{i}\bm{\psi}_{i}(t)
\end{equation}
The independent fundamental modal solutions $\bm{\psi}_{i}(t)$ form a complete functional basis for all possible close-proximity relative motion. The modal constants $\bm{c}$ are functions of the initial conditions $\bm{x}_{0}$, but more generally they can be treated as state variables. 

The modal decomposition is traditionally defined for autonomous linear dynamic systems, but all relative motion except the CW problem is characterized by non-autonomous differential equations. However, the Lyapunov-Floquet (LF) transformation \cite{Nayfeh:1979uq} can be used to equate a linear time-varying (LTV) dynamic system with periodic plant matrix $[A(t)] = [A(t+T)]$ to to a linear time-invariant (LTI) counterpart via a periodic coordinate transformation:
\begin{equation}
\label{LFintro1}
\bm{x} = [P(t)]\bm{z} = [P(t+T)]\bm{z}
\end{equation}
where $\bm{z}$ represents the coordinate set for the LTI equivalent of the system in $\bm{x}$, with the following simple LTI dynamics:
\begin{equation}
\label{LFintro2}
\dot{\bm{z}} = [\Lambda]\bm{z}
\end{equation}
The LF transformation and the LTI matrix are any pair of matrices $[P(t)]$, $[\Lambda]$ satisfying the following matrix differential equation:
\begin{equation}
\label{LFintro3}
[P(t)]^{-1}\left([A(t)][P(t)] - [\dot{P}(t)]\right) = [\Lambda]
\end{equation}

A unique definition of the LF transformation is given below using the monodromy matrix and the state transition matrix (STM). This transformation conveniently equals identity at the epoch time:
\begin{equation}
\label{LFintro4}
[P(t)] = [\Phi(t,t_{0})]e^{-[\Lambda](t - t_{0})}
\end{equation}
\begin{equation}
\label{LFintro5}
[P(t_{0})] = [P(t_{0} + kT)] = [I]
\end{equation}
\begin{equation}
\label{LFintro6}
[\Lambda] = \frac{1}{T}\ln{(\Phi(t_{0}+T,t_{0}))}
\end{equation}
The fundamental modal solutions in this work will thus be of the following form, assuming a diagonalizable LTI plant matrix $[\Lambda]$ with eigenvectors $\bm{v}_{i}$ and eigenvalues $\lambda_{i}$:
\begin{equation}
\label{FS_rm2}
\bm{x}(t) = \sum_{i = 1}^{6}c_{i}[P(t)]\bm{v}_{i}e^{\lambda_{i}t}
\end{equation}
For the special case of purely Keplerian dynamics, $[\Lambda]$ has a Jordan form, so the form of the modal decomposition is slightly modified from above. From the condition that the plant matrix $[A(t)]$ be periodic, it is observed that this perspective applies for relative motion in the vicinity of any periodic chief orbit, regardless of the governing dynamics (two-body, three-body).

It is worth noting that there exists a simple linear mapping of the LF transformations and the resulting LTI matrices when two sets of coordinates are related by a linearized transformation. This relationship is demonstrated below for local coordinates $\bm{x}$ and orbit element differences $\doe$, approximately related by the linearized transformation $[G(t)]$:
\begin{equation}
\label{doeToCart1}
\bm{x} = [G(t)]\doe
\end{equation}
\begin{equation}
\label{transF3}
[P_{\bm{x}}(t)] = [G(t)][P_{\doe}(t)][G(t_{0})]^{-1}
\end{equation}
\begin{equation}
\label{doeToCart4}
[\Lambda_{\bm{x}}] =  [G(t_{0})][\Lambda_{\doe}][G(t_{0})]^{-1}
\end{equation}
See Reference~\citenum{schaub} for the linearized transformation $[G(t)]$.

While there are many possible representations for satellite relative motion, this work mainly uses the popular local Cartesian coordinate representation. The relative state is augmented relative position and velocity $\bm{x} = \left(\bm{\rho}^{\top}, \ \bm{\rho}'^{\top}\right)^{\top}$  resolved in a chief-centered rotating coordinate frame called the Hill or local vertical-local horizontal (LVLH) frame as below:
\begin{subequations}
\label{relpos1}
\begin{align}
\bm{\rho} = & \  x\hat{\bm{e}}_{r} + y\hat{\bm{e}}_{t} + z\hat{\bm{e}}_{n} \\
\bm{\rho}' = & \ \frac{^{H}\text{d}}{\text{d}t}\left(\bm{\rho}\right) = \dot{x}\hat{\bm{e}}_{r} + \dot{y}\hat{\bm{e}}_{t} + \dot{z}\hat{\bm{e}}_{n}
\end{align}
\end{subequations}
The vector triad $\left\{ \hat{\bm{e}}_{r}, \ \hat{\bm{e}}_{t}, \ \hat{\bm{e}}_{n}\right\}$ forming the LVLH frame is defined below in terms of the chief inertial position, velocity, and orbit angular momentum vectors $\bm{r}_{c}$, $\bm{v}_{c}$, and $\bm{h}_{c}$, and $(\ )'$ denotes the time derivative of quantities as seen in this frame.
\begin{subequations}
\label{relpos1}
\begin{align}
\hat{\bm{e}}_{r} = & \ \bm{r}_{c}/r_{c} \\
\hat{\bm{e}}_{n} = & \ \bm{h}_{c}/h_{c} \\
\hat{\bm{e}}_{t} = & \ \hat{\bm{e}}_{n} \times\hat{\bm{e}}_{r} 
\end{align}
\end{subequations}
The general dynamics in local Cartesian coordinates are given in terms of the chief radial vector $\bm{r}$, the angular velocity vector $\bm{\omega}$, and their derivatives in Reference ~\citenum{Casotto}.

To facilitate the description of perturbed and controlled relative motion in terms of the Keplerian relative motion modes, the Keplerian modal decompositions are introduced. There are two of interest. The first is the straightforward solution to the LTI CW equations, and the second is a set of solutions to the general eccentric orbit case $0 < e < 1$. Starting with the CW system, the modal decomposition is computed as in Reference~\citenum{Burnett_AA2021}:
\begin{equation}
\label{xModal_cw}
\bm{x}_{\text{CW}}(t) = \sum_{i=1}^{6}c_{i}\bm{\psi}_{\text{CW}}(t)
\end{equation}
\begin{equation}
\label{Psi_CW}
[\Psi_{\text{CW}}(t)] = \begin{bmatrix} 0 & -\frac{2}{3n} & -\frac{1}{n}\cos{nt} & \frac{1}{n}\sin{nt} & 0 & 0 \\ 
1 & t & \frac{2}{n}\sin{nt} & \frac{2}{n}\cos{nt} & 0 & 0 \\
0 & 0 & 0 & 0 & \frac{2}{n}\sin{nt} & \frac{2}{n}\cos{nt} \\
0 & 0 & \sin{nt} & \cos{nt} & 0 & 0 \\
0 & 1 & 2\cos{nt} & -2\sin{nt} & 0 & 0 \\
0 & 0 & 0 & 0 & 2\cos{nt} & -2\sin{nt} \end{bmatrix}
\end{equation}
\begin{equation}
\label{cvec_CW}
\bm{c} = \begin{pmatrix} y_{0} - \frac{2}{n}\dot{x}_{0} \\ -6nx_{0} - 3\dot{y}_{0} \\ \ 3nx_{0} + 2\dot{y}_{0} \\ \dot{x}_{0} \\ \frac{\dot{z}_{0}}{2} \\ \frac{n}{2}z_{0} \end{pmatrix}
\end{equation}
Because the CW system is already LTI, the LF transformation is just identity. Note that the typical modal decomposition of the CW plant matrix $A = [V][J][V]^{-1}$ would produce complex $[V]$ and thus complex $\bm{c}$. This has been avoided for the derivation of Eq. \eqref{Psi_CW}, and can generally always be avoided. Consider a pair of modes $\bm{\psi}_{m,n} = [P(t)]e^{\lambda_{m,n}t}\bm{v}_{m,n}$ with complex-conjugate eigenvectors $\bm{v}_{m}$, $\bm{v}_{n}$. Taking the constants associated with the complex-conjugate modes $c_{m}, \ c_{n}$ and defining new purely real constants $\tilde{c}_{m} = (c_{m} + c_{n})/2$ and $\tilde{c}_{n} = (c_{m} - c_{n})/2i$, their associated fundamental solutions will have a purely real form.

The preferred decomposition used in this paper for the general Keplerian case (any eccentricity $0<e<1$) is the spherical coordinate modal decomposition, linearly mapped to Cartesian coordinates. This uses the argument of latitude $\theta = \omega+f$ as the independent variable. 
The fundamental modal solutions consist of an along-track rectilinear motion, an offset circular mode, and only two geometrically non-trivial planar modes: a ``teardrop" mode (changing from an initially elliptical shape to a highly pointed shape as eccentricity is increased from zero to one) and a drift mode. The Cartesian modal decomposition, by contrast, has three non-trivial planar modes. There are also two one-dimensional decoupled out-of-plane oscillatory modes. The spherical coordinate modal decomposition is given by the equations below \cite{Burnett_AA2021}, with the $\bm{v}_{i}$ as columns of Eq. \eqref{VforRmatSph}:
\begin{equation}
\label{xModalSol_Msph}
\bm{x}_{c(s)}(\theta) = \sum_{i=1}^{5}c_{i}[F(\theta)]^{-1}[P_{\bm{x}_{s}}(\theta)]\bm{v}_{i} + c_{6}[F(\theta)]^{-1}[P_{\bm{x}_{s}}(\theta)]\left(\bm{v}_{5}(\theta-\theta_{0}) + \bm{v}_{6}\right)
\end{equation}
\begin{equation}
\label{LF_map2}
[P_{\bm{x}_{s}}(\theta)] = [G_{s}(\theta)][P_{\doe}(\theta)][G_{s}(\theta_{0})]^{-1}
\end{equation}
\begin{equation}
\label{PofTheta_QNS1}
[P_{\doe}(\theta)] =  \begin{bmatrix} 1 & 0 & 0 & 0 & 0 & 0 \\ P_{21}(\theta) & P_{22}(\theta) & 0 & P_{24}(\theta) & P_{25}(\theta) & 0 \\ 0 & 0 & 1 & 0 & 0 & 0 \\ 0 & 0 & 0 & 1 & 0 & 0 \\ 0 & 0 & 0 & 0 & 1 & 0 \\ 0 & 0 & 0 & 0 & 0 & 1 \end{bmatrix}
\end{equation}
\begin{subequations}
\label{Pfinal2}
\allowdisplaybreaks
\begin{align}
\begin{split}
P_{21}(\theta) = & \  \frac{\kappa^{2}}{2a}\left(\mathbb{F}_{21}(\theta_{0}) - \mathbb{F}_{21}(\theta)\right) 
\end{split} \\
\begin{split}
\mathbb{F}_{21}(\theta) = & \ \frac{6}{\eta^{3}}\left(\tan^{-1}{\left(\frac{q_{2} + (1 - q_{1})\tan{\left(\frac{\theta}{2}\right)}}{\sqrt{1 - q_{1}^{2} - q_{2}^{2}}}\right)} -\frac{\theta}{2}\right) \\ + & \ \frac{3\left(q_{2} + (q_{1}^{2} + q_{2}^{2})\sin{\theta}\right)}{q_{1}(q_{1}^{2} + q_{2}^{2} - 1)\kappa}
\end{split} \\
\begin{split}
P_{22}(\theta) = & \ \frac{\kappa^{2}}{\kappa_{0}^{2}}
\end{split} \\ 
\begin{split}
P_{24}(\theta) = & \ \frac{\kappa^{2}}{4\left(q_{1}^{2} + q_{2}^{2} - 1\right)}\left(\mathbb{F}_{24}(\theta_{0}) - \mathbb{F}_{24}(\theta)\right)
\end{split} \\
\begin{split}
\mathbb{F}_{24}(\theta) = & \ \frac{4(q_{2} + \sin{\theta})}{\kappa^{2}} + \frac{4\sin{\theta}}{\kappa}
\end{split} \\
\begin{split}
P_{25}(\theta) = & \ \frac{\kappa^{2}}{4\left(q_{1}^{2} + q_{2}^{2} - 1\right)}\left(\mathbb{F}_{25}(\theta_{0}) - \mathbb{F}_{25}(\theta)\right) 
\end{split} \\
\begin{split}
\mathbb{F}_{25}(\theta) = & \ \frac{4(1 - q_{1}^{2} + q_{2}\sin{\theta})}{q_{1}\kappa^{2}}  + \frac{4q_{2}\sin{\theta}}{q_{1}\kappa}
\end{split} 
\end{align}
\end{subequations}
\begin{equation}
\label{FofTheta}
[F(\theta)] = 
\begin{bmatrix} 
1 & 0 & 0 & 0 & 0 & 0  \\
0 & 1/r & 0 & 0 & 0 & 0 \\
0 & 0 & 1/r & 0 & 0 & 0 \\
0 & 0 & 0 & 1 & 0 & 0 \\
0 & -\dot{r}/r^{2} & 0 & 0 & 1/r & 0  \\ 
0 & 0 & -\dot{r}/r^{2} & 0 & 0 & 1/r
\end{bmatrix}
\end{equation}
\begin{equation}
\label{GsofTheta}
[G_{s}(\theta)] = 
\begin{bmatrix} \frac{r}{a} & \frac{v_{r}}{v_{t}}r & 0 & -\frac{r}{p}(2aq_{1} + r\cos{\theta}) & -\frac{r}{p}(2aq_{2} + r\sin{\theta}) & 0 \\ 0 & 1 & 0 & 0 & 0 & \cos{i} \\ 0 & 0 & \sin{\theta} & 0 & 0 & -\cos{\theta}\sin{i} \\ -\frac{v_{r}}{2a} & \left(\frac{1}{r} - \frac{1}{p}\right)h & 0 & \frac{1}{p}(v_{r}aq_{1} + h\sin{\theta}) & \frac{1}{p}(v_{r}aq_{2} - h\cos{\theta}) & 0 \\ -\frac{3\dot{\theta}}{2a} & -2\frac{v_{r}}{r} & 0 & \frac{\dot{\theta}}{p}(3aq_{1} + 2r\cos{\theta}) & \frac{\dot{\theta}}{p}(3aq_{2} + 2r\sin{\theta}) & 0 \\ 0 & 0 & \dot{\theta}\cos{\theta} & 0 & 0 & \dot{\theta}\sin{\theta}\sin{i} \end{bmatrix}
\end{equation}
\begin{equation}
\label{VforRmatSph}
[V_{R_{\bm{x}_{s}}}] = \begin{bmatrix} 0 & 0 & 0 & 0 & \frac{2R_{21}a}{\gamma}AC\gamma a & 0 \\
1 & 0 & 0 & 0 & \frac{2R_{21}a}{\gamma}(B+1)^{2}C & 0 \\
0 & 1 & 0 & 0 & 0 & 0 \\
0 & 0 & 1 & 0 & \frac{2R_{21}a}{\gamma}B\gamma a & 0 \\
0 & 0 & -\frac{A}{\gamma a} & 0 & -\frac{4R_{21}a}{\gamma}A(B+1) & 1 \\
0 & 0 & 0 & 1 & 0 & 0
\end{bmatrix}
\end{equation}
In the above equations, $v_{r} = \dot{r}$ and $v_{t} = r\dot{\theta}$, $q_{1} = e\cos{\omega}$, $q_{2} = e\sin{\omega}$, and the shorthand $\text{s}$ and $\text{c}$ are sin and cos. 
Remaining quantities are defined below:
\begin{subequations}
\allowdisplaybreaks
\label{xcQuant}
\begin{align}
R_{21} = & \ -\frac{3a\eta}{2r_{0}^{2}}\\
A = & \ q_{2}\cos{\theta_{0}} - q_{1}\sin{\theta_{0}} \\
B = & \ q_{1}\cos{\theta_{0}} + q_{2}\sin{\theta_{0}}\\
C = & \ \frac{hr_{0}^{2}}{a\mu\gamma} \\
\gamma = & \ A^{2} + B^{2} - 1 \\
\eta = & \ \sqrt{1 - q_{1}^{2} - q_{2}^{2}} \\
\kappa = & \ 1 + q_{1}\cos{\theta} + q_{2}\sin{\theta} \\
\kappa_{0} = & \ 1 + q_{1}\cos{\theta_{0}} + q_{2}\sin{\theta_{0}} \\
\end{align}
\end{subequations}
The fundamental solution constants, which are quite important, are given in terms of initial spherical relative coordinates by Eq. \eqref{c123sph}:
\begin{subequations}
\label{c123sph}
\begin{align}
c_{1} = & \ -\frac{v_{t,0}}{v_{r,0}r}\delta r_{0} + \theta_{r,0} \\
c_{2} = & \ \phi_{r,0} \\
c_{3} = & \ \frac{1}{C}\left(\frac{\left(1 - \frac{r_{0}}{p}\right)v_{t,0}}{v_{r,0}}\delta r_{0} + C\delta \dot{r}_{0} \right) \\
c_{4} = & \ \dot{\phi}_{r,0} \\
c_{5} = & \  -\frac{v_{t,0}}{3v_{r,0}a}n\left(\frac{r_{0}}{p}\right)\delta r_{0} \\
c_{6} = & \ \frac{\mu}{hr_{0}^{2}}\left(1 + \frac{p}{r_{0}}\right)\delta r_{0} + \frac{v_{r,0}}{v_{t,0}r_{0}}\delta\dot{r}_{0} + \dot{\theta}_{r,0}
\end{align}
\end{subequations}
The function $[F(\theta)]$ linearly maps the Cartesian coordinates $x, \ y, \ z, \ \dot{x}, \ \dot{y}, \ \dot{z}$ to the spherical relative state coordinates $\delta r, \theta_{r}, \phi_{r}, \delta\dot{r}, \dot{\theta}_{r}, \dot{\phi}_{r}$, and can be used for getting the spherical coordinate initial conditions required by Eq.~\eqref{c123sph}.
Note that singularities exist with the given Keplerian modal formulation, appearing when $q_{1} = e\cos{\omega} = 0$ or when  $e\sin{f_{0}} = 0$. These can be avoided by selecting nearby conditions, setting the offending terms to a small number $\epsilon$ instead of exactly zero.

The modal decompositions can always be normalized to maximize ease in geometric interpretation of the relative scales of the $c_{i}$s:
\begin{equation}
\label{xRescale}
\bm{x}(t) = \sum_{i=1}^{6}c_{i}\bm{\psi}_{i}(t) = \sum_{i=1}^{6}c_{i}\|\bm{\psi}_{i}\|\hat{\bm{\psi}}_{i}(t) = \sum_{i=1}^{6}\overline{c}_{i}\hat{\bm{\psi}}_{i}(t) 
\end{equation}
where the choice of normalization is left to the reader. A simple choice is to make it so that the maximum relative range of a normalized mode over the course of one orbit is unity. This is the standard normalization scheme applied for all results in this paper.

\begin{figure}[h!]
\centering
\includegraphics[]{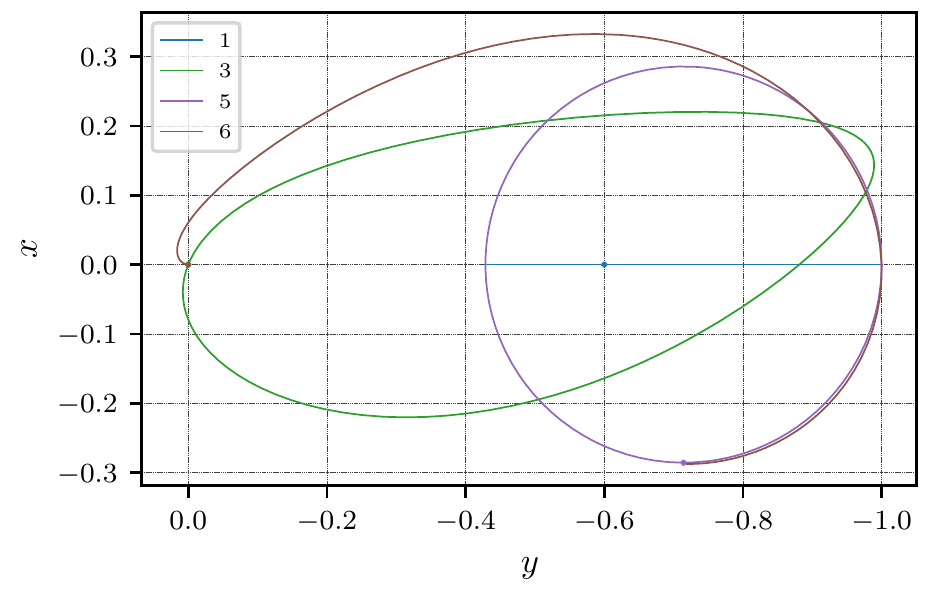}
\caption{Example In-plane Normalized Relative Motion Modes} \vspace{-5mm}
\label{fig:modeSphEx1}
\end{figure}
\begin{figure}[h!]
\centering
\includegraphics[]{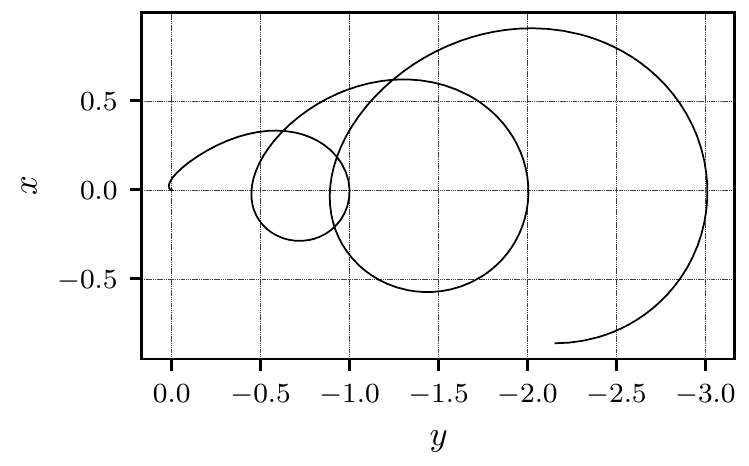}
\caption{Drift Mode, 3 Orbits} \vspace{-5mm}
\label{fig:driftModeSphEx1}
\end{figure}

The normalized planar ($x$ and $y$) modes for the linearly mapped spherical coordinates are produced in Figure \ref{fig:modeSphEx1} for an example case with an unperturbed chief orbit about Earth with elements $a = 12000$ km, $e = 0.4$, $i = 25^{\circ}$, $\Omega = 0^{\circ}$, $\omega = 270.001^{\circ}$, $f_{0} = 90^{\circ}$. The initial points for each mode are indicated with a dot. The out-of-plane modes, which are simply linearly independent oscillations in $z$, are not shown. From Figure \ref{fig:modeSphEx1}, the planar projection of all close-proximity unperturbed relative motion in this case is revealed to be a linear sum of geometrically simple independent motions. First is an along-track oscillation (mode 1, shown flipped about $y = 0$ for convenience), a distorted ``teardrop" mode (mode 3), an offset circle mode (mode 5), and a drift mode (mode 6). Note that the chief eccentricity dictates the ratio of the offset circle's radius to its distance from the origin, and the along-track oscillation of mode 1 has the same maximum and minimum $y$ bounds as the offset circle mode. Note also that the drift mode behaves as a composition of a secularly expanding mode 5 and a translation along the $y$ axis -- see the last term in Eq. \eqref{xModalSol_Msph}. To better illustrate this behavior, the drift mode is shown for 3 chief orbit periods in Figure \ref{fig:driftModeSphEx1}. 
While the uncontrolled state evolves as a sum of the modes, the state can be uniquely parameterized by the stationary $\bm{c}$ and the time since epoch $t - t_{0}$. In this sense, the modal solution constants function similarly to orbit elements, but with the exact linear mapping $\bm{x}(t) = [\Psi(t)]\bm{c}$ to local relative state coordinates. There are many possible alternate choices of fundamental solutions $\bm{\psi}_{i}(t)$, but Lyapunov-Floquet theory is typically useful for identifying geometrically favorable solutions.
It is worth noting that the use of the modal decomposition approach is not limited to Keplerian dynamics. This method can be applied to study motion in the vicinity of any periodic orbit (and work is ongoing to extend this to the more general case of a quasi-periodic chief orbit). This will be shown with an application of the method to the relative motion problem in the circular restricted three-body problem.

\subsection{The Modal Constants as State Variables}
For the Keplerian problem, any uncontrolled close-proximity relative motion $\bm{x}(t)$ can be expressed as a constant weighted sum of the fundamental modal motions. In the event that the orbital dynamics are perturbed, or control is active, the general behavior can still be represented by a time-varying weighted sum of the modal motions:
\begin{subequations}
\label{VOP1}
\begin{align}
\dot{\bm{x}} = & \ \bm{f}(\bm{x},\bm{u},t) \\
\bm{x}(t) = & \ [\Psi(t)]\bm{c}(t)
\end{align}
\end{subequations}
where the first equation gives the general true nonlinear dynamics, and $\bm{c}(t)$ is made to vary such that these dynamics are satisfied.
The vector of constants is allowed to vary in time as $\bm{c}(t)$ such that any motion $\bm{x}(t)$ can be represented. This is enabled by the following osculating condition:
\begin{subequations}
\label{VOP2}
\begin{align}
\dot{\bm{x}}(t) = & \ \frac{\partial\bm{x}}{\partial t} + \frac{\partial\bm{x}}{\partial\bm{c}}\dot{\bm{c}} = \bm{f}(\bm{x},\bm{u},t) \\
\frac{\partial\bm{x}}{\partial t} = & \ [A(\oec^{*},t)]\bm{x}
\end{align}
\end{subequations}
where $[A(\oec^{*},t)]$ is the plant matrix for the original linearized dynamics used to compute the fundamental modal solutions $[\Psi(t)]$.
The following dynamics satisfy the osculating condition:
\begin{equation}
\label{VOP3}
\dot{\bm{c}} = [\Psi(t)]^{-1}\left(\bm{f}(\bm{x},\bm{u},t) - [A(\oec^{*},t)]\bm{x}\right)
\end{equation}
Assuming that $\bm{x}(t)$ remains in the linear regime, Eq. \eqref{VOP3} reduces to linear dynamics in $\bm{c}$:
\begin{equation}
\label{VOP3}
\dot{\bm{c}} = [\Psi(t)]^{-1}\Big([A(t)] - [A^{*}(t)]\Big)[\Psi(t)]\bm{c} + [\Psi(t)]^{-1}[B_{\bm{x}}]\bm{u}
\end{equation}
where $[B_{\bm{x}}] = [0_{3\times3} \ I_{3\times3}]^{\top}$ if $\bm{x}$ is in Cartesian coordinates, and $[A(t)]$ is the plant matrix for linearization about the true (deviated) chief orbit. In the case that the orbital dynamics are unchanged, Eq. \eqref{VOP3} obtains a form where $\bm{c}$ is only influenced by control:
\begin{equation}
\label{VOP3_alt}
\dot{\bm{c}} = [\Psi(t)]^{-1}[B_{\bm{x}}]\bm{u}
\end{equation}
Because $\bm{c}$ has no linearized plant matrix, it is possible in this case to design control to track a desired natural reference trajectory $\bm{x}_{r}(t)$ using a regulation controller in $\bm{c}$ space, because the dynamics of the error $\delta\bm{c} = \bm{c} - \bm{c}_{r}$ are of the same form as Eq. \eqref{VOP3_alt}.

In the case that the true orbit is weakly perturbed in comparison to the orbit used to develop the modal decomposition, one may write the true plant matrix in terms of the nominal plus a small deviation:
\begin{equation}
\label{VOP4}
[A(t)] \approx [A^{*}(t)] + \varepsilon[\delta\tilde{A}(t)]
\end{equation}
where $|\varepsilon| \ll 1$ is a small parameter. In this case, Eq. \eqref{VOP3} is rewritten:
\begin{equation}
\label{VOP5}
\begin{split}
\dot{\bm{c}} = & \ \varepsilon[\Psi(t)]^{-1}[\delta\tilde{A}(t)][\Psi(t)]\bm{c} + [\Psi(t)]^{-1}[B_{\bm{x}}]\bm{u} \\ \equiv & \ \varepsilon[\Omega(t)]\bm{c} + [B_{\bm{c}}(t)]\bm{u}
\end{split}
\end{equation}
Eq. \eqref{VOP5} is interesting because as the relative state is written as $\bm{x}(t) = [\Psi(t)]\bm{c}(t)$, the perturbed close-proximity dynamics of relative motion can be factored into a traditional slow-fast dynamical system:
\begin{subequations}
\label{SF1}
\begin{align}
[\dot{\Psi}] = & \ [A^{*}(t)]\Psi \\
\dot{\bm{c}} = & \ \varepsilon [\Omega(t)]\bm{c} + [B_{\bm{c}}(t)]\bm{u}
\end{align}
\end{subequations}
where from Eq. \eqref{VOP5}, the matrices $[\Omega]$ and $[B_{\bm{c}}]$ are also functions of the modal solutions. 
The ``slowness" of the dynamics in $\bm{c}$ implies that for weakly perturbed cases, the state $\bm{c}(t)$ will change slowly from $\bm{c}(t_{0})$ in comparison to how the $\bm{\psi}_{i}(t)$ change from $\bm{\psi}_{i}(t_{0})$ and how $\bm{x}(t)$ changes from $\bm{x}(t_{0})$. This allows interpretation of the evolving relative motion using the osculating modal constants, similarly to how the osculating orbital elements are used to study general perturbed satellite orbits.

\subsection{Extension to the Circular Restricted Three-Body Problem}
Eqs.~\eqref{FS_rm1} - \eqref{FS_rm2} apply for relative motion in the vicinity of any periodic orbit, not just in the Keplerian case. Thus, this paper explores relative motion in the circular restricted three-body problem (CR3BP) from the same perspective. For this problem, the states and time are non-dimensionalized in the traditional manner \cite{Koon:2006rf}, with dynamics resolved in the rotating coordinates in which the two primaries appear fixed. The relative motion in the vicinity of a periodic orbit is studied in CR3BP rotating coordinates and is given as $\delta\bm{x} = [P(\tau)]\bm{z}$ (where the delta differentiates this from the usual LVLH frame relative state), and $\dot{\bm{z}} = [\Lambda]\bm{z}$. 

Due to the Hamiltonian nature of the CR3BP, the monodromy matrix has a repeated eigenvalue of $\lambda_{m,1} = \lambda_{m,2} = +1$. Furthermore, it can be shown that the state rate $\dot{\bm{X}}(\tau_{0}) = \dot{\bm{X}}(\tau_{0}+T)$ is an eigenvector corresponding to the unity eigenvalue. There is not a second eigenvector -- the unity eigenvalue has algebraic multiplicity two and geometric multiplicity one. Solving for the generalized eigenvector $\bm{v}_{2}$ is straightforward:
\begin{equation}
\label{GE3bp1}
([\Phi(\tau_{0}+T,\tau_{0})] - [I_{6\times6}])\bm{v}_{2} = \bm{v}_{1}
\end{equation}
Then, the resulting Jordan decomposition of the monodromy matrix is as below:
\begin{equation}
\label{MonoJord1}
[\Phi(\tau_{0}+T,\tau_{0})] = [V_{m}][J_{m}][V_{m}]^{-1} = [\bm{v}_{1} \ \bm{v}_{2} \ldots \bm{v}_{6}]\begin{bmatrix} 1 & 1 & 0 & \ldots \\ 0 & 1 & 0 & \ldots \\ 0 & 0 & \lambda_{m,3} \\ \vdots & \vdots & & \ddots \end{bmatrix}[\bm{v}_{1} \ \bm{v}_{2} \ldots \bm{v}_{6}]^{-1}
\end{equation}
The LTI matrix is given by $[\Lambda] = \frac{1}{T}\ln(\Phi(\tau_{0}+T,\tau_{0}))$:
\begin{equation}
\label{MonoJord2}
[\Lambda] = \frac{1}{T}[V_{m}]\ln(J_{m})[V_{m}]^{-1} = \frac{1}{T}[\bm{v}_{1} \ \bm{v}_{2} \ldots \bm{v}_{6}]\begin{bmatrix} 0 & 1 & 0 & \ldots \\ 0 & 0 & 0 & \ldots \\ 0 & 0 & \ln(\lambda_{m,3}) \\ \vdots & \vdots & & \ddots \end{bmatrix}[\bm{v}_{1} \ \bm{v}_{2} \ldots \bm{v}_{6}]^{-1}
\end{equation}
The solution to linearized relative motion in the CR3BP has the following form, with a drift mode as a result of the defectiveness:
\begin{equation}
\label{CR3BP_dxsol1}
\delta\bm{x}(\tau) = c_{1}[P(\tau)]\bm{v}_{1} + c_{2}[P(\tau)]\left(\bm{v}_{1}\tau + \bm{v}_{2}\right) + c_{3}\bm{\psi}_{3}(\tau) + \ldots + c_{6}\bm{\psi}_{6}(\tau)
\end{equation}
where the trivial mode is the first listed mode, being periodic in CR3BP coordinates, and the drift mode is listed second. Then there are four other modes (stable, unstable, or center) starting with $\bm{\psi}_{3}$. Recall that the delta notation is adopted for the CR3BP case to differentiate these differential CR3BP coordinates from the LVLH coordinates in the Keplerian case.

The defectiveness of the LTI form, its double-zero eigenvalues, and the resulting secular drift mode might remind the reader of the Keplerian relative motion modal decomposition. The drift mode has a simple physical interpretation: bounded relative motion must satisfy a period-matching condition, so motion on nearby orbits of different periods would violate this condition and result in nonzero projection into a local drift mode in the linearized system. 
This relationship also appears in the Keplerian relative motion problem: bounded purely along-track relative motion is possible (representing points of a different phase along the orbit), and in the case that the no-drift condition is violated, the drift occurs along this same along-track direction. 

The above discussion of the defectiveness of the monodromy matrix (and the underlying LTI form) and the resulting drift mode enables a fully analytic view of the relative motion modes. 
First, consider the case of relative motion in the vicinity of a stable orbit, where all the relative motion modes are bounded (i.e. all eigenvalues of $[\Lambda]$ have zero real part). In this case, there is one pair of trivial eigenvalues and there are two non-trivial pairs of eigenvalues:
\begin{subequations}
\label{Eig_ex1}
\begin{align}
\lambda_{1,2} = & \ 0 \\
\lambda_{3,4} = & \ \pm i\omega_{1} \\
\lambda_{5,6} = & \ \pm i\omega_{2}
\end{align}
\end{subequations}
Let the complex-conjugate eigenvector pairs associated with the frequencies $\omega_{1}$ and $\omega_{2}$ be written as $\bm{v}_{3,4} = \bm{v}_{\mathbb{R}_{1}} \pm i\bm{v}_{\mathbb{I}_{1}}$ and $\bm{v}_{5,6} = \bm{v}_{\mathbb{R}_{2}} \pm i\bm{v}_{\mathbb{I}_{2}}$, respectively. Re-factoring, the following form is obtained for the modal decomposition of the relative motion:
\begin{subequations}
\label{PsiEx1p1}
\begin{align}
\delta\bm{x}(\tau) = & \ \sum_{i=1}^{6}c_{i}\bm{\psi}_{i}(\tau) \\
\bm{\psi}_{i} = & \ [P(\tau)]\bm{\eta}_{i}(\tau)
\end{align}
\end{subequations}
\begin{subequations}
\label{eta1T6ex1}
\begin{align}
\bm{\eta}_{1}(\tau) = & \  \bm{v}_{1}\\
\bm{\eta}_{2}(\tau) = & \  \bm{v}_{1}\tau + \bm{v}_{2} \\
\bm{\eta}_{3}(\tau) = & \ 2\left(\bm{v}_{\mathbb{R}_{1}}\cos{\left(\omega_{1}\tau\right)} - \bm{v}_{\mathbb{I}_{1}}\sin{\left(\omega_{1}\tau\right)}\right) \\
\bm{\eta}_{4}(\tau) = & \ -2\left(\bm{v}_{\mathbb{R}_{1}}\sin{\left(\omega_{1}\tau\right)} + \bm{v}_{\mathbb{I}_{1}}\cos{\left(\omega_{1}\tau\right)}\right) \\
\bm{\eta}_{5}(\tau) = & \ 2\left(\bm{v}_{\mathbb{R}_{2}}\cos{\left(\omega_{2}\tau\right)} - \bm{v}_{\mathbb{I}_{2}}\sin{\left(\omega_{2}\tau\right)}\right) \\
\bm{\eta}_{6}(\tau) = & \ -2\left(\bm{v}_{\mathbb{R}_{2}}\sin{\left(\omega_{2}\tau\right)} + \bm{v}_{\mathbb{I}_{2}}\cos{\left(\omega_{2}\tau\right)}\right)
\end{align}
\end{subequations}
where the $\bm{c}$ is given as a function of $\delta\bm{x}_{0}$:
\begin{equation}
\label{cfromx0_ex1}
\bm{c} = [\overline{V}]^{-1}\delta\bm{x}_{0}
\end{equation}
\begin{equation}
\label{Vbar_ex1}
[\overline{V}] = [\bm{v}_{1}, \ \bm{v}_{2}, \ 2\bm{v}_{\mathbb{R}_{1}}, \ -2\bm{v}_{\mathbb{I}_{1}}, \ 2\bm{v}_{\mathbb{R}_{2}}, \ -2\bm{v}_{\mathbb{I}_{2}}]
\end{equation}

Because modes $\bm{\psi}_{3}$ - $\bm{\psi}_{6}$ are generally composed of multiple incommensurate frequencies, they trace out complex and unintuitive shapes on long timespans. The trivial modes associated with the double-zero eigenvalues of $[\Lambda]$ are comparatively simple, because they are $T$-periodic. Together, these modes form the basis of all close-proximity relative motion in the vicinity of the periodic orbit. 

Another example is the case of two trivial modes, two center modes, a stable mode, and an unstable mode. The modal decomposition is given as below, where $\bm{v}_{3,4} = \bm{v}_{\mathbb{R}_{1}} \pm i\bm{v}_{\mathbb{I}_{1}}$:
\begin{subequations}
\label{PsiEx2p1}
\begin{align}
\delta\bm{x}(\tau) = & \ \sum_{i=1}^{6}c_{i}\bm{\psi}_{i}(\tau) \\
\bm{\psi}_{i} = & \ [P(\tau)]\bm{\eta}_{i}(\tau)
\end{align}
\end{subequations}
\begin{subequations}
\label{eta1T6ex2}
\begin{align}
\bm{\eta}_{1}(\tau) = & \  \bm{v}_{1} \\
\bm{\eta}_{2}(\tau) = & \  \bm{v}_{1}\tau + \bm{v}_{2} \\
\bm{\eta}_{3}(\tau) = & \ 2\left(\bm{v}_{\mathbb{R}_{1}}\cos{\left(\omega_{1}\tau\right)} - \bm{v}_{\mathbb{I}_{1}}\sin{\left(\omega_{1}\tau\right)}\right) \\
\bm{\eta}_{4}(\tau) = & \ -2\left(\bm{v}_{\mathbb{R}_{1}}\sin{\left(\omega_{1}\tau\right)} + \bm{v}_{\mathbb{I}_{1}}\cos{\left(\omega_{1}\tau\right)}\right) \\
\bm{\eta}_{5}(\tau) = & \ \bm{v}_{5}e^{\lambda_{5}\tau} \\
\bm{\eta}_{6}(\tau) = & \ \bm{v}_{5}e^{\lambda_{6}\tau}
\end{align}
\end{subequations}
\begin{equation}
\label{Vbar_ex2}
[\overline{V}] = [\bm{v}_{1}, \ \bm{v}_{2}, \ 2\bm{v}_{\mathbb{R}_{1}}, \ -2\bm{v}_{\mathbb{I}_{1}}, \ \bm{v}_{5}, \ \bm{v}_{6}]
\end{equation}
In this case, the existence of an unstable mode generally results in relative motion being unstable if there is any projection of $\delta\bm{x}$ into the unstable subspace. All other dynamic cases can be explored as needed in the same manner as the above two cases.

\section{Applications of the Fundamental Solution Constants}
\subsection{Analytic Benefits as a State Representation}
For researchers interested in approximating satellite relative motion efficiently, it is worth noting that first-order perturbative expansions of the relative motion in $\bm{c}$ space are quite convenient in comparison to the typical investigations in $\bm{x}$ space or in $\doe$ coordinates:
\begin{equation}
\label{Pertdyn4}
\dot{\bm{c}} =  \varepsilon [\Omega(t)]\bm{c}
\end{equation}
\begin{equation}
\label{Pertdyn5}
\bm{c}(t) \approx \bm{c}^{0} + \varepsilon\bm{c}^{1}(t)
\end{equation}
\begin{equation}
\label{Pertdyn6}
\dot{\bm{c}}^{1} = [\Omega(t)]\bm{c}^{0}
\end{equation}
where $\bm{c}^{1}(t)$ is solved simply by integration, and $\bm{x}$ is expanded in terms of the perturbed $\bm{c}$:
\begin{equation}
\label{Pertdyn7}
\bm{c}^{1}(t) = \int_{0}^{t}[\Omega(\varphi)]\text{d}\varphi \ \bm{c}^{0}
\end{equation}
\begin{equation}
\label{Pertdyn8}
\bm{x}(t) \approx [\Psi(t)]\left([I_{6\times6}] + \varepsilon\int_{0}^{t}[\Omega(\varphi)]\text{d}\varphi\right)[\Psi(0)]^{-1}\bm{x}_{0}
\end{equation}

\subsection{Continuous Control Using the Fundamental Solution Constants}

In the classical Linear Quadratic Tracking (LQT) problem, the formulation in terms of local coordinates, such as the Cartesian representation, is inconvenient. To execute tracking control, the $6\times6$ gain matrix $[K]$ must first be propagated backwards via a Riccati matrix differential equation (requiring simultaneous back-propagation of the chief orbit and computation of the plant matrix $[A]$), and the 6-dimensional co-state $\bm{s}(t)$ must also be back-propagated. 

Consider instead that the fundamental solution constants $\bm{c}$ are used as the state representation, and the desired trajectory to track is natural, thus $\bm{u}_{r} = \bm{0}$ and $\bm{c}_{r}(t) = \bm{c}_{r}(t_{0})$. The state error is $\delta\bm{c} = \bm{c} - \bm{c}_{r}$ with the following simple dynamics, assuming that the chief orbit is unperturbed:
\begin{equation}
\label{LQR_c1}
\delta\dot{\bm{c}} = [B_{\bm{c}}(t)]\bm{u}
\end{equation}
Furthermore, the LQT problem in $\bm{x}$ space reduces to the Linear Quadratic Regulator (LQR) problem in $\bm{c}$ space, with cost function, optimal control, and simplified Riccati equation below:
\begin{equation}
\label{LQT_c2}
\tilde{J}_{\text{LQR}} =   \frac{1}{2}\delta\bm{c}(t_{f})^{\top}[S]\delta\bm{c}(t_{f}) + \frac{1}{2}\int_{t_{0}}^{t_{f}}\Big(\delta\bm{c}(t)^{\top}[Q]\delta\bm{c}(t) + \bm{u}^{\top}[R]\bm{u} \Big)\text{d}t
\end{equation}
\begin{equation}
\label{LQT_c3}
\bm{u}(t) = -[R]^{-1}[B_{\bm{c}}]^{\top}[K]\delta\bm{c}(t)
\end{equation}
\begin{equation}
\label{LQT_c4}
[\dot{K}] = [K][B_{\bm{c}}(t)][R]^{-1}[B_{\bm{c}}(t)]^{\top}[K] - [Q], \ \ [K(t_{f})] = [S]
\end{equation}
With this formulation, the controlled $\bm{x}(t)$ will track natural trajectory $\bm{x}_{r}(t)$ through control in $\bm{c}$ space -- where there is no need to back-propagate any co-state dynamics, and the Riccati equation is also greatly simplified by the absence of an $[A]$ matrix. However, a complication is that the choice of satisfactory gains is not as straightforward in $\bm{c}$ space as it is in $\bm{x}$ space.

\subsection{Impulsive Control Using the Fundamental Solution Constants}
The relative motion parameterization in terms of fundamental solution constants is well-suited for impulsive maneuver-based control strategies. Returning to Eq. \eqref{SF1}, in the absence of disturbances, $[\Omega(t)] = [0_{6\times6}]$, and the solution for $\bm{c}$ can be expressed in terms of a series of impulsive maneuvers $\Delta\bm{v}_{i} = \Delta\bm{v}(t_{i})$:
\begin{equation}
\label{SF1_disc_up}
\bm{c} = \bm{c}_{0} + \sum_{i=1}^{N}[B_{c}(t_{i})]\Delta\bm{v}_{i}
\end{equation}
Consider the optimal control problem of minimizing the total delta-V subject to the dynamics in Eq. \eqref{SF1_disc_up}:
\begin{equation}
\label{JfromU}
J = \sum_{i=1}^{N}\Delta v_{i}
\end{equation}
where $\Delta v_{i} = \|\Delta\bm{v}(t_{i})\|$. Also, let $\Delta\bm{c} = \bm{c} - \bm{c}_{0}$. Some useful results are borrowed from Reference~\citenum{Guffanti_ICs}. First, consider the sets of control inputs $\bm{u}(t) \in \mathscr{U}$ and the reachable variations $\Delta\bm{c} \in \mathscr{C}$ with cost no greater than $J$:
\begin{equation}
\label{Uset1}
\mathscr{U}(J) = \left\{ \bm{u}(t) : \bm{u}(t) = \sum_{i=1}^{N}\Delta\bm{v}_{i}, \sum_{i=1}^{N}\Delta v_{i} \leq J \right\}
\end{equation}
\begin{equation}
\label{Cset1}
\mathscr{C}(J) = \left\{ \Delta\bm{c} : \Delta\bm{c} = \sum_{i=1}^{N}[B_{c}(t_{i})]\Delta\bm{v}_{i}, \sum_{i=1}^{N}\Delta v_{i} \leq J \right\}
\end{equation}
As discussed in Reference~\citenum{Guffanti_ICs}, the set $\mathscr{C}$ is compact and convex, and scales linearly with $J$. Furthermore, for a minimum cost $J_{\text{min}}$ to achieve a desired variation, the desired difference in constants $\Delta\bm{c}^{*} = \bm{c}^{*} - \bm{c}_{0}$ lies on the boundary of the set. The minimum delta-V to reach this goal in $N$ maneuvers can be obtained in terms of the unit vector $\hat{\bm{\eta}}$, which is normal to the boundary of $\mathscr{C}$ at $\Delta\bm{c}^{*}$:
\begin{equation}
\label{Jmin1}
J_{\text{min}} = \frac{\hat{\bm{\eta}}^{\top}\Delta\bm{c}^{*}}{\text{max}_{t_{i}\in[t_{0},t_{f}]}\|\hat{\bm{\eta}}^{\top}[B_{c}(t_{i})]\|}
\end{equation}
For integration constant control formulations, Reference~\citenum{Guffanti_ICs} described a means of numerically obtaining $\hat{\bm{\eta}}$ using a convex solver, then linearly solving for an optimal sequence of $N \leq n$ impulsive maneuvers for a dynamic system with $n$ state variables. In general, for the formation flying problem, a minimum of two maneuvers are required. For the unperturbed problem, only control action induces movement in $\mathscr{C}$ -- the flow of the integrable dynamics has no effect. This property, combined with the compactness and convexity of $\mathscr{C}$, allows for powerful geometric interpretations for the fuel-optimal impulsive maneuver problem. However, any significant perturbations will play a disruptive role, inducing drifts in $\Delta\bm{c}$ that would need to be compensated. 

The algorithm for solving for an optimal maneuver sequence is given below for the unperturbed problem:\cite{Guffanti_ICs}
\begin{enumerate}
\allowdisplaybreaks
\item Solve the following second-order cone program for the optimal value $\bm{\eta}^{*}$:
\begin{equation}
\label{Guffanti_J}
\begin{split}
& \text{maximize} \ \ \ \ \tilde{J} = \bm{\eta}^{\top}\Delta\bm{c}^{*} \\
& \text{subject to} \ \ \ \ \|[B_{c}](t)^{\top}\bm{\eta}\| \leq 1 \ \text{for} \ t \in [t_{0}, \ldots, t_{j}, \ldots, t_{f}]
\end{split}
\end{equation}
where $[t_{0}, \ldots, t_{j}, \ldots, t_{f}] \in \mathbb{R}^{k}_{\geq 0}$ is a chosen discretization of the control interval. 
\item Determine all times $t_{i} \in [t_{0}, \ldots, t_{j}, \ldots, t_{f}]$ for which $|\|[B_{c}(t)^{\top}\bm{\eta}^{*}\| - 1| < \epsilon$ for some tolerance $\epsilon \ll 1$. This will yield an $N$-maneuver sequence, with $N  \ll k$, for which the $i$\textsuperscript{th} impulse is directed along the unit vector:
\begin{equation}
\label{impulse_i}
\hat{\bm{u}}_{i} = [B_{c}](t)^{\top}\bm{\eta}^{*}
\end{equation}
\item The set of delta-V maneuvers $\left\{\Delta v_{i}\right\}$ needs to satisfy the linear system of equations:
\begin{equation}
\label{SF1_disc_up}
\sum_{i=1}^{N}[B_{c}(t_{i})]\hat{\bm{u}}_{i}\cdot\Delta v_{i} = \Delta\bm{c}^{*}
\end{equation}
\end{enumerate}

A traditional relative motion control approach is to use a factorization of the STM to uniquely solve for a 2-burn delta-V sequence given some specified initial and final times. This could also be used, with the desired target relative orbit efficiently identified in $\bm{c}$ space. 
This $\bm{c}^{*}$ would then be mapped to $\bm{x}^{*}$ via the fundamental solutions. 
By necessity, a two-burn solution using either approach will obtain the same answer.

\section{Numerical Results for the Keplerian Case}
\subsection{Keplerian Relative Orbit and its Modes}
\begin{figure}[h!]
\centering
\includegraphics[]{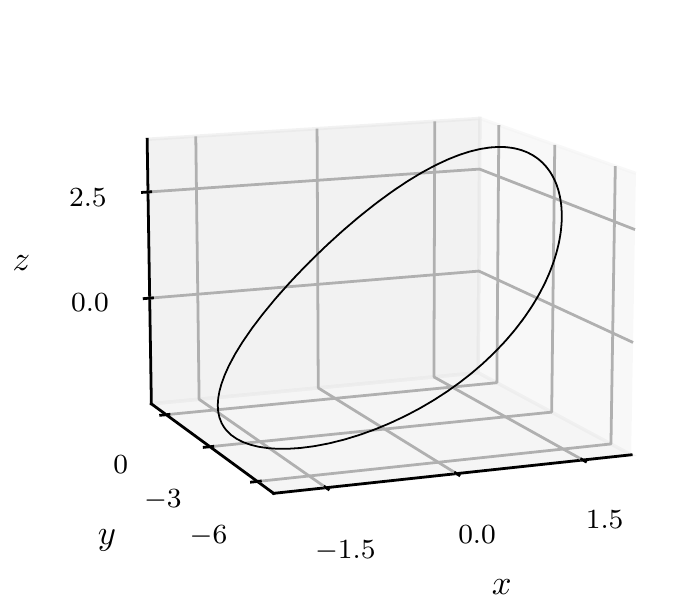}
\caption{Example Relative Orbit} \vspace{-5mm}
\label{fig:RO_3d}
\end{figure}

\begin{figure}[h!]
\centering
\includegraphics[]{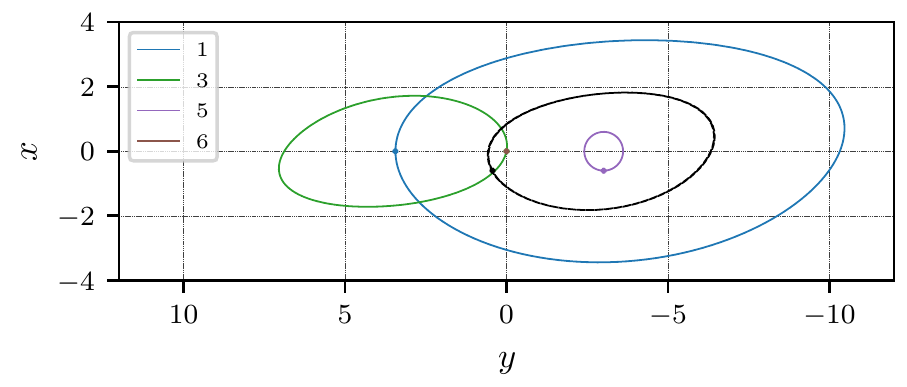}
\caption{Cartesian Planar Modes} \vspace{-5mm}
\label{fig:Modes2Dcart}
\end{figure}

\begin{figure}[h!]
\centering
\includegraphics[]{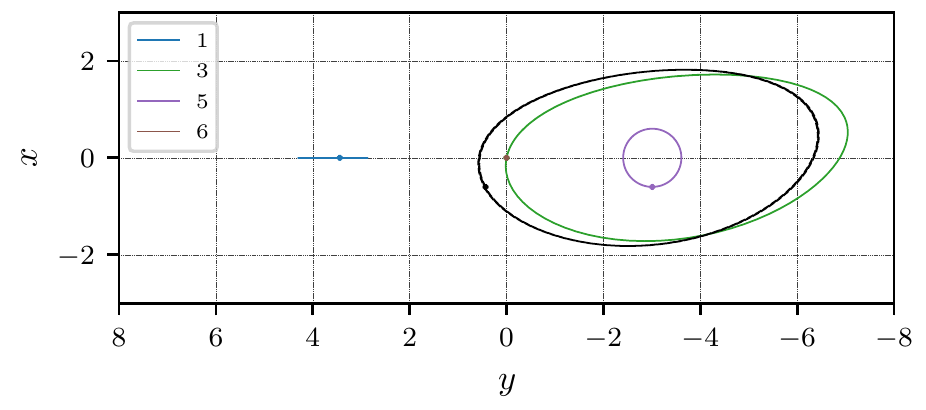}
\caption{Spherical Planar Modes} \vspace{-5mm}
\label{fig:Modes2Dsph}
\end{figure}

\begin{figure}[h!]
\centering
\includegraphics[]{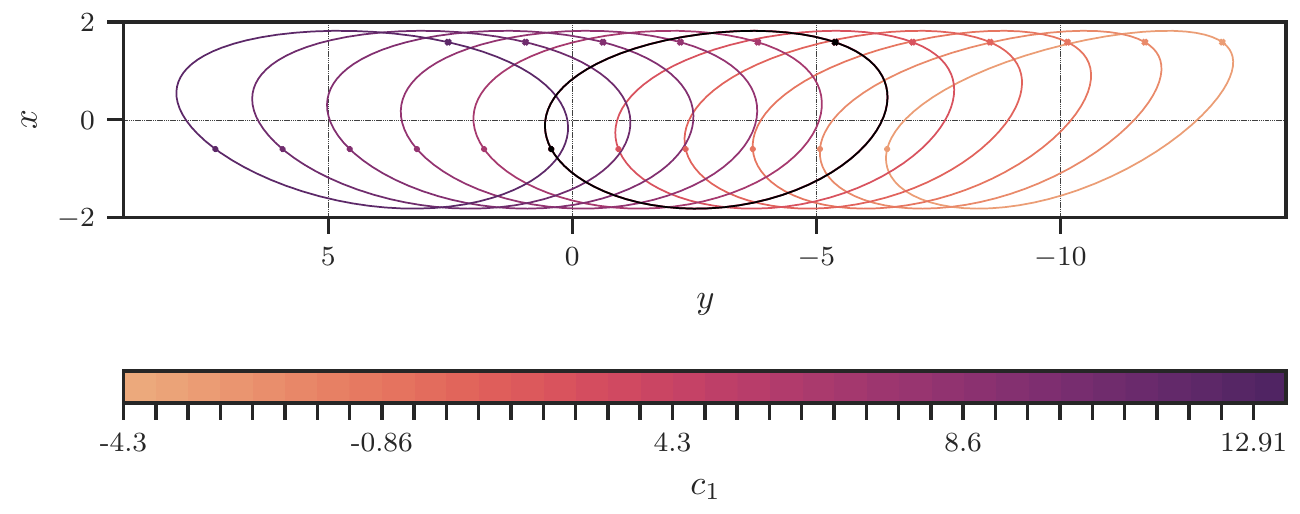}
\caption{Variations of the Planar Relative Motion with $c_{1}$} \vspace{-5mm}
\label{fig:c1varRO2d}
\end{figure}

Consider the example of a bounded relative orbit in the vicinity of an eccentric chief orbit, given by the data listed in Table \ref{table:RO_params}. Note that while $\delta a = 0$, the drift constant $c_{6}$ does not exactly equal zero. This is because $\delta a = 0$ is the nonlinear no-drift requirement, and it is not perfectly captured by the linearized no-drift condition. This is a well-known property of linearized solutions -- consider for example that the linearized no-drift condition for the CW system is $\dot{y} + 2nx = 0$, which only linearly approximates $\delta a = 0$. From the data in Table \ref{table:RO_params}, the resulting Keplerian relative orbit is depicted in 3D in Figure \ref{fig:RO_3d}, and the planar projection appears as a black closed curve in Figures \ref{fig:Modes2Dcart} -- \ref{fig:c1varRO2d}.  
\begin{table}[h!]
\centering
\caption{Simulation Parameters}
\begin{tabular}{l|l}
Parameter             & Value                                                    \\ \hline
Chief Orbit Elements  & $\oec = \left(a, e, i, \Omega, \omega, f_{0} \right) =  \left(8600\text{ km}, 0.2, 25^{\circ}, 0^{\circ}, 270.001^{\circ}, 90^{\circ}\right)$ \\
Deputy Relative Orbit & $\doe = \oec_{d} - \oec_{c} = \left(0.0, 0.0002, 0.02^{\circ}, 0^{\circ}, 0^{\circ}, 0.003^{\circ} \right)$ \\
Modal Constants & $\bm{c} = \left(4.3, 0.0, 7.07, 3.60, 3.61, -0.014\right)$         
\end{tabular} \vspace{-5mm}
\label{table:RO_params}
\end{table}

Figure \ref{fig:Modes2Dcart} shows the modal decomposition of the in-plane component of the relative motion using the modes developed in Cartesian coordinates. Figure \ref{fig:Modes2Dsph} shows the modal decomposition of the in-plane motion using the modes developed in spherical coordinates and linearly mapped to Cartesian coordinates. For both plots, the initial point of the orbiter and the initial point in each mode are given by dots. Because the relative orbit is non-drifting, the drift mode contribution is zero -- thus the mode appears as a non-moving point at the origin. For both modal decompositions, the modes shown sum linearly to reproduce the observed relative motion in black. In other words, $\bm{x}_{\text{2D}}(t) = c_{1}\bm{\psi}_{\text{2D},1} + c_{3}\bm{\psi}_{\text{2D},3} + c_{5}\bm{\psi}_{\text{2D},5}$. Recall that the out-of-plane modes (2 and 4) have no in-plane component -- they exist only in $z$, and are completely decoupled from the in-plane modes.

Comparing Figures \ref{fig:Modes2Dcart} and \ref{fig:Modes2Dsph}, the spherical coordinate-based modal decomposition reproduces the true relative orbit in a much more straightforward and intuitive manner than the Cartesian coordinate-based counterpart. The motion is represented as a sum of a rectilinear along-track motion (mode 1), a distorted elliptical motion (mode 3), and the offset circular trajectory (mode 5). This is the simplest geometric representation possible for Keplerian relative motion, with two of the three bounded in-plane motions given as basic shapes. 

\begin{figure}[h!]
\centering
\includegraphics[]{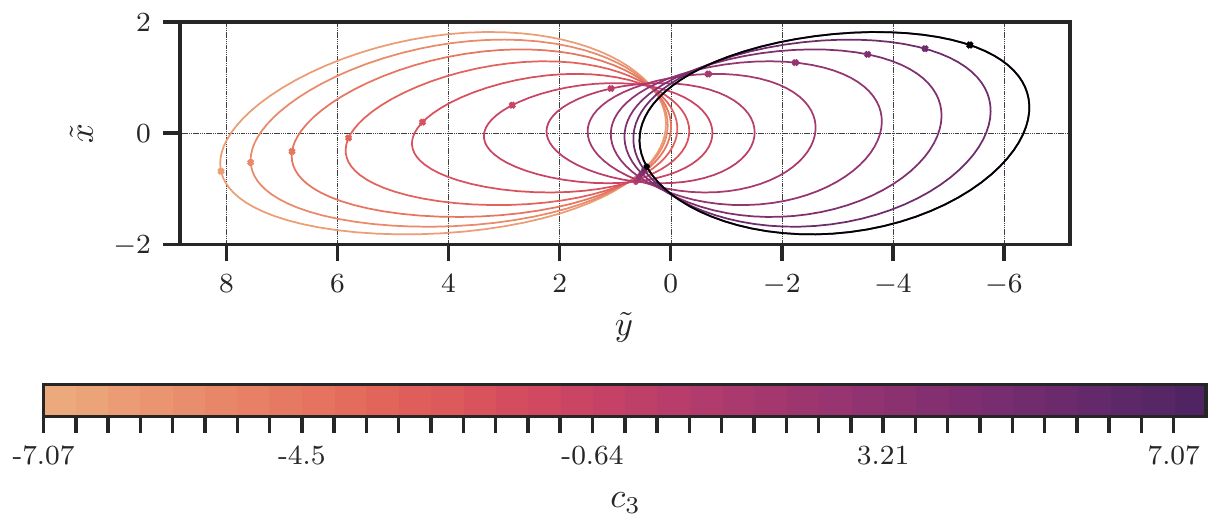}
\caption{Variations of the Planar Relative Motion with $c_{3}$ (Re-scaled)} \vspace{-5mm}
\label{fig:c3varRO2d}
\end{figure}
\begin{figure}[h!]
\centering
\includegraphics[]{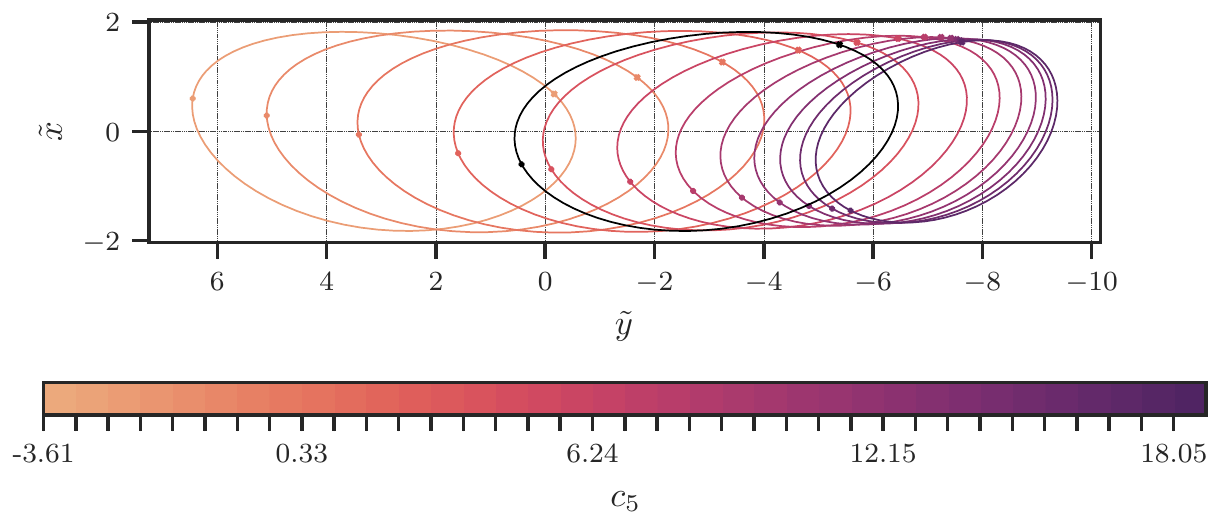}
\caption{Variations of the Planar Relative Motion with $c_{5}$ (Re-scaled)} \vspace{-7mm}
\label{fig:c5varRO2d}
\end{figure}
Because the fundamental modal motions only need to be computed once, variations in the relative motion due to changes in the modal constants can be explored with a minimal amount of numerical computation. For example, exploring a range of variations in $c_{1}$, Figure \ref{fig:c1varRO2d} is produced. The initial ($t = t_{0}$) and mid-orbit ($t = t_{0} + \tfrac{T}{2}$) points are denoted with a ``$\bullet$" and with an ``\texttt{x}", respectively, and the original relative orbit is given in black. The effect of isolated changes in $c_{1}$ is to shift the motion further along in the along-track direction as $c_{1}$ is increased, with the additional effect of rotating and distorting the planar component of the relative orbit. Note that the $x$ scale of the relative orbit is not affected at all. Similar figures can be generated to isolate the effects of changes in $c_{3}$ and $c_{5}$ on the relative orbit shape and location. However, changing the scales of $c_{3}$ and $c_{5}$ also change the size of the relative orbit. To display the characteristic changes in relative orbits with these parameters clearly on individual plots, the relative orbits are computed across desired ranges for these parameters, as was done for Figure \ref{fig:c1varRO2d}, then the orbits are re-scaled such that $\|\bm{c}_{\text{new}}\| = \|\bm{c}_{\text{old}}\|$ to preserve the original relative orbit scale. The resulting plots are given in Figure \ref{fig:c3varRO2d} for variations in $c_{3}$ and Figure \ref{fig:c5varRO2d} for variations in $c_{5}$. 

In Figure \ref{fig:c3varRO2d}, as $c_{3}$ is decreased from its original value of 7.07, the re-scaled relative orbit shifts from the original relative orbit (given in black) to more centered and symmetric relative motion in the middle of the range (near $c_{3} = 0$), to an essentially reversed version of the original for $c_{3} < 0$. Note that there would also be accompanying relative orbit scale changes with changing value of $c_{3}$, but the re-scaled orbit plot sacrifices this information to better show the variations in the relative orbit geometry. Figure \ref{fig:c5varRO2d} shows the variations in re-scaled relative orbit due to changes in $c_{5}$, with an original value of $c_{5} = 3.61$. The negative value is essentially flipped about the $x$-axis, and as the value is increased, the re-scaled relative orbits gradually circularize as the contribution of the circular mode 5 is increased in relative scale. 

Figures \ref{fig:c1varRO2d} -- \ref{fig:c5varRO2d} show that the parameter space for the in-plane component of bounded relative orbits is only three-dimensional. The two out-of-plane modes add an additional two dimensions -- completely decoupled from the in-plane design. From the perspective of the modal constants, it is conceptually easy and numerically efficient for the astrodynamicist to explore all possible useful types of relative motion that can exist. In this manner, the vector of modal constants $\bm{c}$ serves as the design space, and also uniquely determines the relative motion state when combined with a given time since epoch $t - t_{0}$. As discussed earlier, it is also possible to compute how the constants vary under the influence of non-Keplerian dynamics. With such a study, the influence of perturbations on relative motion can be viewed as an evolving alteration of the relative scales of the constituent Keplerian relative motion modes that form the basis for the unperturbed problem.

\subsection{Effects of Perturbations -- Modeling with $\mathbf{J_{2}}$}
To demonstrate the behavior of the Keplerian modal constants under the influence of perturbations, consider the ubiquitous example of $J_{2}$-perturbed relative motion, which highly relevant for Earth orbits. To compute the perturbed dynamics of the Keplerian modal constants, Eq. \eqref{VOP3} is used with equations from Reference~\citenum{Casotto} providing the linearized perturbed relative motion dynamics, evaluated using the following equations for the acceleration, differential acceleration, and jerk induced by the $J_{2}$-perturbed Keplerian dynamics: 
\begin{equation}
\label{rddot_J2}
\ddot{\bm{r}} = -\frac{\mu}{r^{3}}\bm{r} - \frac{3\mu J_{2} R^{2}}{2r^{4}}\left(\left(1-5\left(\hat{\bm{r}}\cdot \hat{\mathbf{K}}\right)^{2}\right)\hat{\bm{r}}+ 2\left(\hat{\bm{r}}\cdot \hat{\mathbf{K}}\right)\hat{\mathbf{K}}\right)
\end{equation}
\begin{equation}
\label{delrddot_J2}
\begin{split}
\nabla_{\bm{r}}\ddot{\bm{r}} = & -\frac{3\mu J_{2} R^{2}}{2r^{5}}\bigg[ \left( 1 - 5\left( \hat{\mathbf{K}}\cdot \hat{\mathbf{r}} \right)^{2}\right)\text{I} + 2\hat{\mathbf{K}} \hat{\mathbf{K}}^{\top} + 5 \left( 7\left( \hat{\mathbf{K}}\cdot \hat{\mathbf{r}} \right)^{2} - 1\right)\hat{\mathbf{r}}\hat{\mathbf{r}}^{\top}\\
& \quad - 10 \left( \hat{\mathbf{K}} \cdot \hat{\mathbf{r}} \right) \left( \hat{\mathbf{K}}\hat{\mathbf{r}}^{\top} + \hat{\mathbf{r}}\hat{\mathbf{K}}^{\top} \right) \bigg]
\end{split}
\end{equation}
\begin{equation}
\label{rdddot_J2}
\dddot{\bm{r}} = \left(\nabla_{\bm{r}}\ddot{\bm{r}}\right)\dot{\bm{r}}
\end{equation}
where $\hat{\bm{K}}$ denotes the polar axis unit vector, $R$ is the equatorial radius, $\text{I}$ is the $3\times 3$ identity matrix, and $\hat{\bm{r}} = \bm{r}/r$. The dynamics given by Eq. \eqref{VOP3} are integrated in parallel with the $J_{2}$-perturbed chief orbit. 

\begin{figure}[h!]
\centering
\includegraphics[]{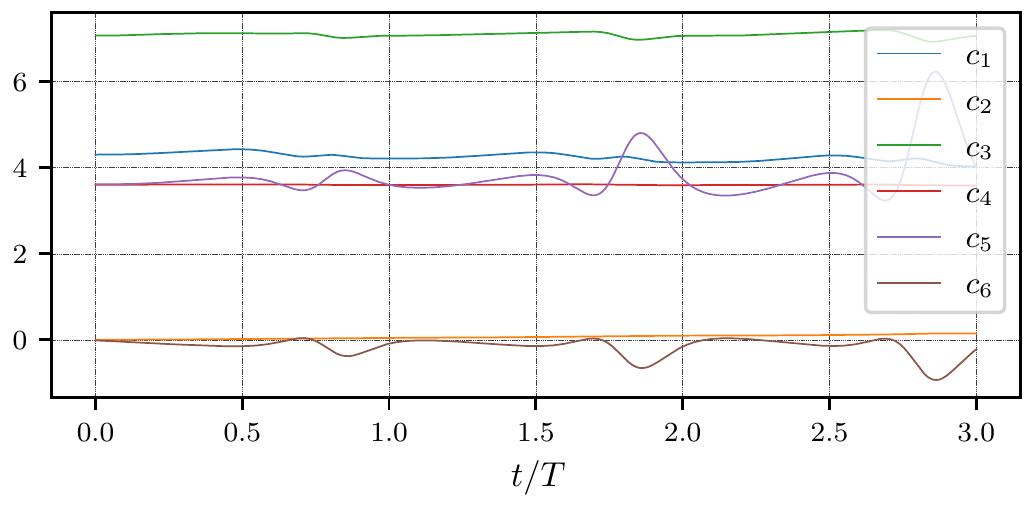}
\caption{Variation of Modal Constants with $J_{2}$} \vspace{-5mm}
\label{fig:cvecJ2_1}
\end{figure}
The same initial chief orbit and deputy relative orbit conditions from Table \ref{table:RO_params} are selected, but with Earth's $J_{2}$ perturbation active. As a result of this perturbation, variations are induced in the Keplerian modal constants. Figure \ref{fig:cvecJ2_1} shows the resulting behavior in $\bm{c}$ for 3 unperturbed chief orbit periods. The effect of $J_{2}$ is limited to small oscillations in the modal constants, but these oscillations grow over time, which is an unfortunate but unavoidable property. For the case of $J_{2}$, it seems that the modified orbital frequency due to the perturbation requires that the secular drift mode (mode 6) be used to fully describe the perturbed state. This is because all other modes are periodic on the interval $[0, T]$, and variations in their sums would be unable to describe a relative orbit on the shortened interval $[0, T - \Delta T_{J_{2}}]$. Additionally, the $J_{2}$ perturbation induces slow long-term drift in the relative orbit. The drift mode is used to describe the perturbed solution, and it grows and shifts over time, so variations in other modes (primarily mode 5) are induced to compensate for these variations. This yields the opposing behaviors of $c_{5}$ and $c_{6}$ seen in Figure \ref{fig:cvecJ2_1}. Despite these growing oscillations, the long-term drift in $c_{5}$ and $c_{6}$ and in the other $c_{i}$ parameters is quite slow. The averaged dynamics in $\bm{c}$ could thus be a useful lens for studying relative motion in the perturbed problem, especially for even zonal harmonics like $J_{2}$, but such a study is not explored here. Lastly, the growing oscillations are not a major problem for modeling, because the fundamental solutions can always be re-initialized as needed. 

\subsection{Impulsive Maneuver Control Examples -- Keplerian Orbit}
\begin{figure}[h!]
\centering
\includegraphics[]{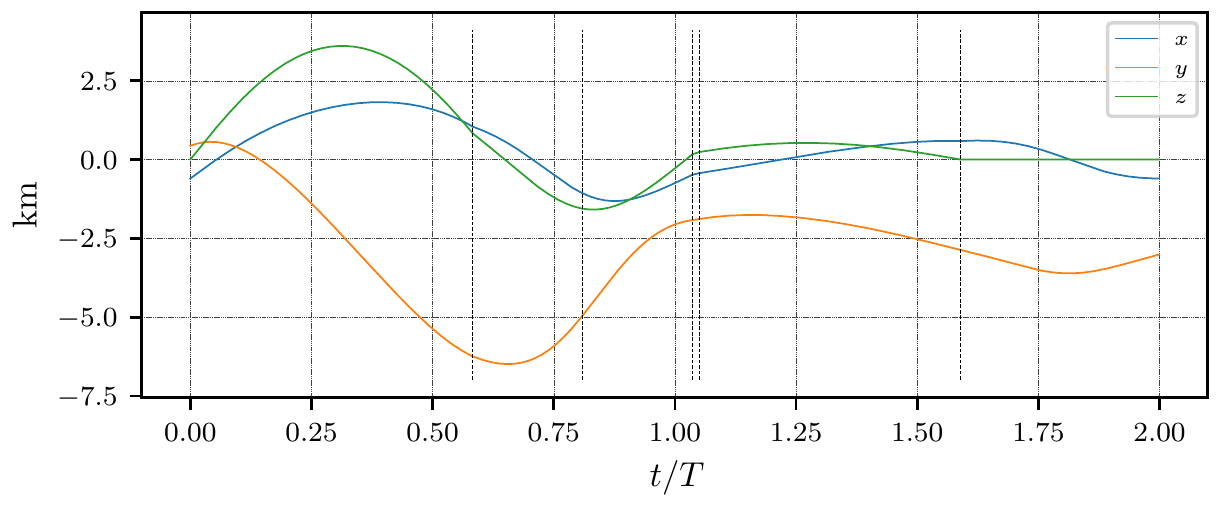}
\caption{Five-Burn Maneuver, Local Cartesian Coordinates} \vspace{-5mm}
\label{fig:xyz_ImpCont2}
\end{figure}

\begin{figure}[h!]
\centering
\includegraphics[]{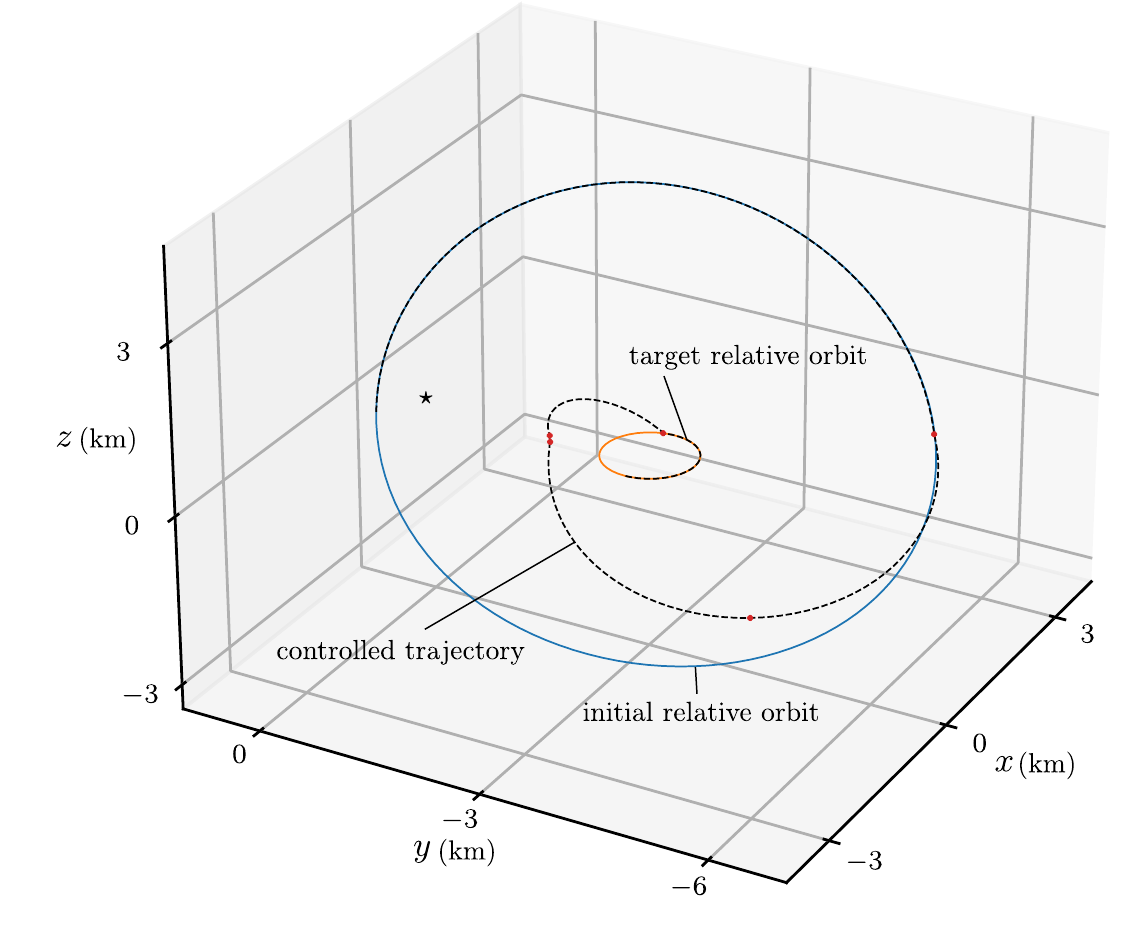}
\caption{Relative Orbit Transfer Using the Modal Constants} 
\label{fig:xyz3D_ImpCont2}
\end{figure}

To demonstrate unperturbed control using the Keplerian modal constants, consider the problem of changing from the initial relative motion dictated by the data in Table \ref{table:RO_params} to a new planar non-drifting relative orbit parameterized by $c_{5} = 3.61$. 
To implement this test, the previously discussed impulsive maneuver-based control solution strategy is implemented in Python using \texttt{cvxpy}. 
\begin{figure}[h!]
\centering
\includegraphics[]{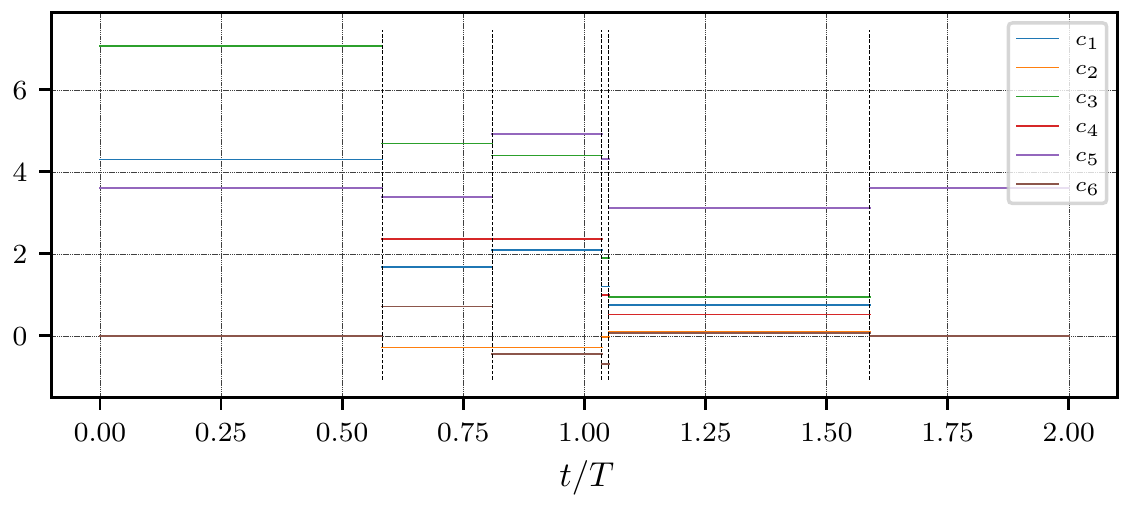}
\caption{Five-Burn Maneuver, Modal Constants} 
\label{fig:cvec_ImpCont2}
\end{figure}
For this impulsive control example with the Keplerian case, time is discretized into 100 points on the interval $t_{i} \in [t_{0} + 1590.6 \ \text{s}, \ t_{0} + 12724.7 \ \text{s}]$. Thus, the number of constraints given by Eq. \eqref{Guffanti_J} is 100. This does not stress the solver, and the optimal maneuver sequence is found quite quickly. The resulting impulsive control solution consists of 5 maneuvers for a combined delta-V of only 2.7 m/s, compared to 7.0 m/s for a two-burn transfer in the same interval. The changes in $\bm{c}$ with each maneuver are plotted in Figure \ref{fig:cvec_ImpCont2}, with vertical dotted lines indicating each maneuver. The relative position components are plotted in Figure \ref{fig:xyz_ImpCont2}. The initial, transfer, and final relative motions are shown in 3D in Figure \ref{fig:xyz3D_ImpCont2}. The initial relative orbit is the large blue closed curve, and the final relative orbit is the small orange planar circular trajectory. The transfer trajectory is given by the dashed black line, and the maneuver points are indicated with red dots. The chief location is indicated by a star. 

Control using the modal solution constants is highly convenient, due to both the straightforward geometric interpretation of the $c_{i}$ parameters and the efficient means by which multi-maneuver impulsive control schedules can be obtained. This is demonstrated by the preceding simple examples with Keplerian dynamics. However, extending this control design to account for the effect of perturbations on $\bm{c}$ is necessary for elegant flight implementation -- particularly for long-duration control maneuver sequences. This will be explored in future work, but mitigating the effects of perturbations in control design in $\bm{c}$ space should be straightforward, because even in the perturbed relative motion case, the $\bm{c}$ parameters do not lose their geometric meaning. This is highly convenient, and generally does not hold for most other perturbed integral representations. For example, for relative motion parameterizations using orbit element differences $\doe = \oec_{d} - \oec_{c}$, the perturbations modify $\oec_{c}$, and as a result, the resulting exact local coordinate behavior mapped from a particular desired $\doe$ changes over time. 
\section{Numerical Examples for the CR3BP Case}
\begin{figure}[h!]
\centering
\includegraphics[]{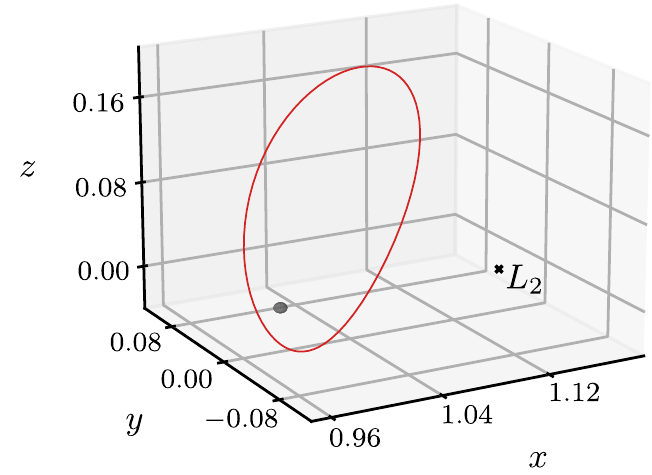}
\caption{Stable $L_{2}$ Halo Orbit} 
\label{fig:L2HaloOrbEx1}
\end{figure}

\begin{figure}[h!]
    \centering
    \subfloat[Mode 3]{{\includegraphics[]{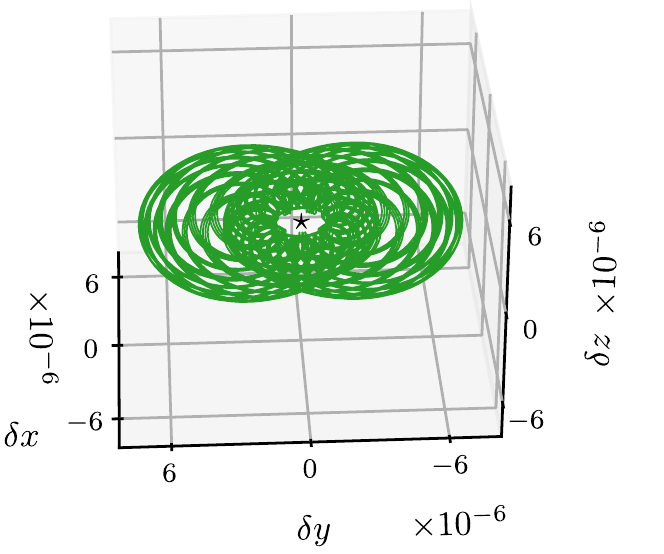} }}
    \subfloat[Mode 5]{{\includegraphics[]{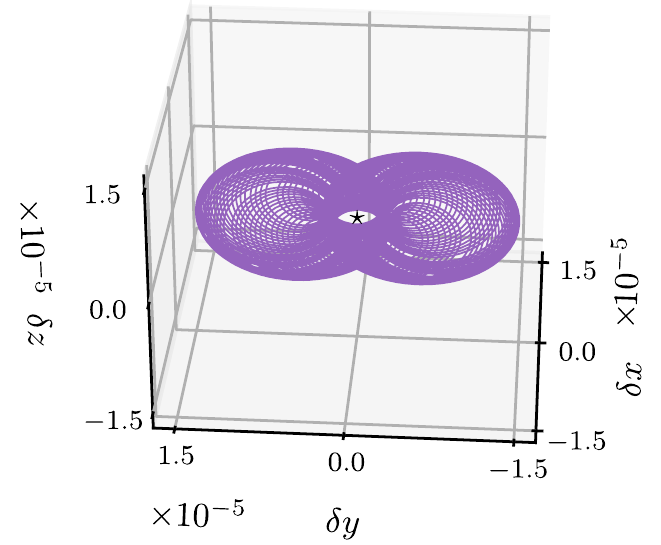} }}%
    \caption{Center Modes, Stable $L_{2}$ Northern Halo Orbit} 
    \label{fig:CenterModesEx1}%
\end{figure}

Consider a stable northern $L_{2}$ halo orbit with an orbit period of $T = 9.504$ days. This orbit is given in Figure~\ref{fig:L2HaloOrbEx1}. For this orbit, there are 4 center modes and two trivial modes ($\bm{\psi}_{1}$ and $\bm{\psi}_{2}$) in its vicinity. The center modes are composed of incommensurate frequencies, so they trace out complex shapes over long timespans. This is demonstrated with plots of $\bm{\psi}_{3}$ and $\bm{\psi}_{5}$ propagated with the linearized dynamics for 240 chief orbits, given by Figure~\ref{fig:CenterModesEx1}. The scale shown corresponds to relative motion on the km scale, but is plotted in the dimensionless CR3BP length scale. The dimensionless frequencies are $\omega_{1} = 1.2511$ and $\omega_{2} = 0.7604$. Not shown is the trivial mode, which traces a closed curve with each chief orbit. Note that the chief location is plotted with a star.

\begin{figure}[h!]
\centering
\includegraphics[]{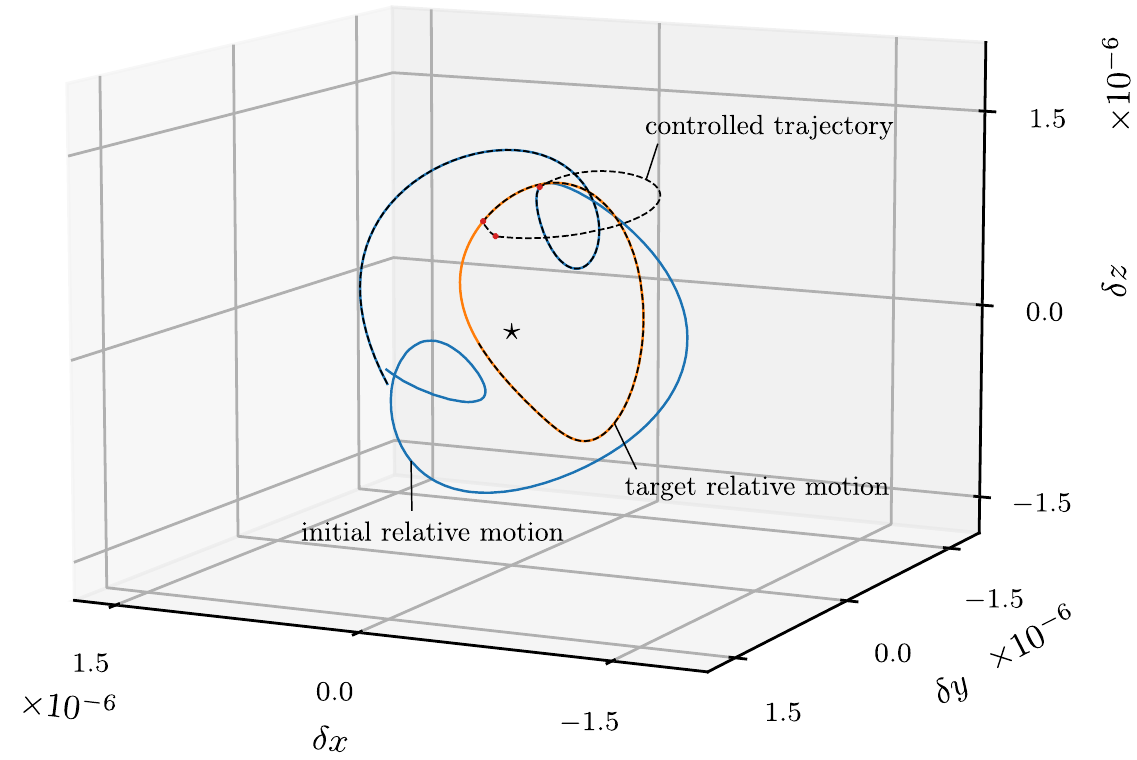}
\caption{Relative Motion with Impulsive Control, Stable $L_{2}$ Halo Orbit} 
\label{fig:HaloControlEx1}
\end{figure}

As a demonstration of the impulsive control strategy discussed earlier, consider the control case summarized in Table~\ref{table:L2HaloEx1}. The initial motion is bounded but irregular, and the bounded trivial mode is targeted. The resulting relative motion is plotted in Figure~\ref{fig:HaloControlEx1}. The uncontrolled trajectory is given in blue for two chief orbits, the target trajectory is in orange, and the controlled trajectory is given by the dashed line, with impulsive maneuver points marked by red dots. The chief is shown as a star. This figure is plotted in the rotating CR3BP coordinates, not the LVLH frame. The figure shows that the impulsive maneuver-based control strategy is successful in achieving the desired relative motion.

\begin{table}[h!]
\centering
\caption{Halo Orbit Control Simulation Parameters (Stable Orbit Example)} 
\label{table:L2HaloEx1}
\begin{tabular}{l|l}
Parameter               & Value                                                                        \\ \hline
Initial relative motion & $\bm{c}_{0} = (0, 0, 0.2, 0.1, 0.08, 0)\alpha$, $\alpha = 5.2\times 10^{-6}$ \\
Initial state           & $\delta\bm{x}_{0} = (-0.01, 0.309, -0.005, 0.168, -0.002, 0.362)\alpha$      \\
Desired relative motion & $\bm{c}^{*} = (0.2, 0, 0, 0, 0, 0)\alpha$                                    \\
Maneuver interval       & $\tau \in [1.23, 3.29]$ \\ 
Maneuver times & $\tau_{1} = 1.37$, $\tau_{2} = 2.876$, $\tau_{3} = 3.013$ \\
Resulting maneuvers     & $\Delta\bm{v}_{1} = (-0.281, 0.094, 0.161)\alpha$ \\ & $\Delta\bm{v}_{2} = (-0.131, -0.027, 0.085)\alpha$ \\ & $\Delta\bm{v}_{3} = (-0.221, -0.064, 0.121)\alpha$ 
\end{tabular} 
\end{table}

\begin{figure}[h!]
    \centering
    \subfloat[Halo Orbit, Example 2(a)]{{\includegraphics[]{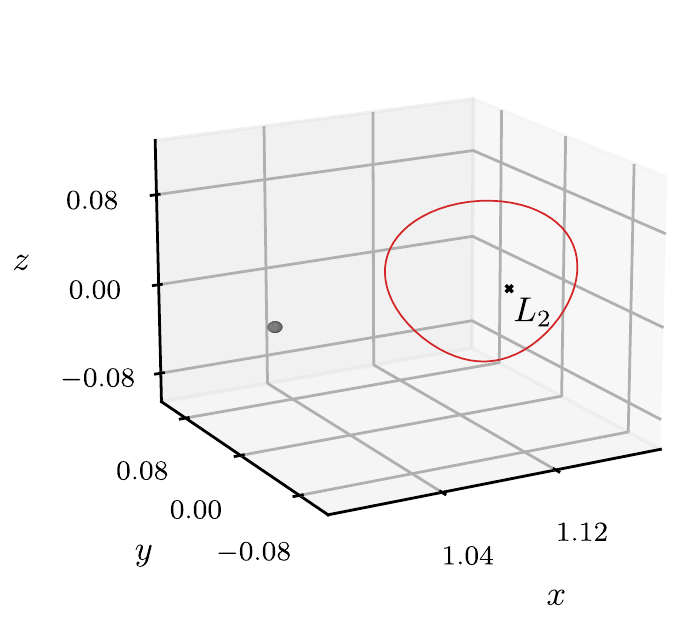} }}
    \subfloat[Center Mode]{{\includegraphics[]{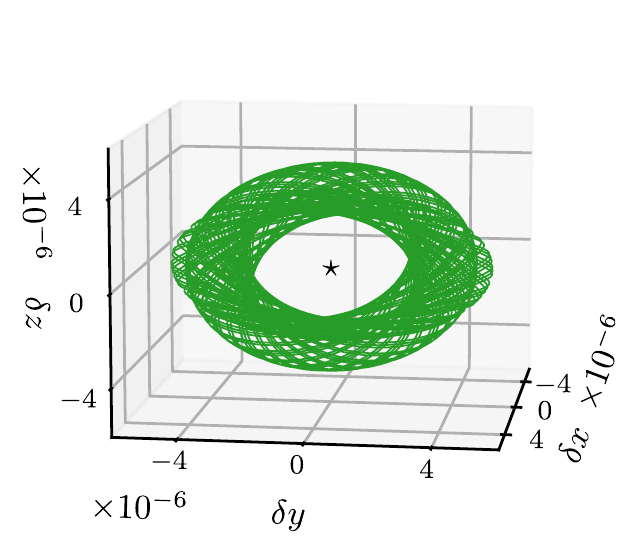} }}%
    \caption{Unstable $L_{2}$ Northern Halo Orbit and a Center Mode}%
    \label{fig:L2Ex2}%
\end{figure}

Next, consider two cases of control near an unstable $L_{2}$ halo orbit. This is a northern $L_{2}$ halo orbit with a period of $T = 14.676$ days. 
For this orbit, there are the trivial bounded and drift modes $\bm{\psi}_{1}$ and $\bm{\psi}_{2}$, two center modes $\bm{\psi}_{3}$ and $\bm{\psi}_{4}$, a stable mode $\bm{\psi}_{5}$, and an unstable mode $\bm{\psi}_{6}$. 
The center mode $\bm{\psi}_{3}$ is propagated for many orbits and given also in Figure~\ref{fig:L2Ex2}.  Its dimensionless frequency is $\omega_{1} = 0.1288$.
\begin{figure}[h!]
    \centering
    \includegraphics[]{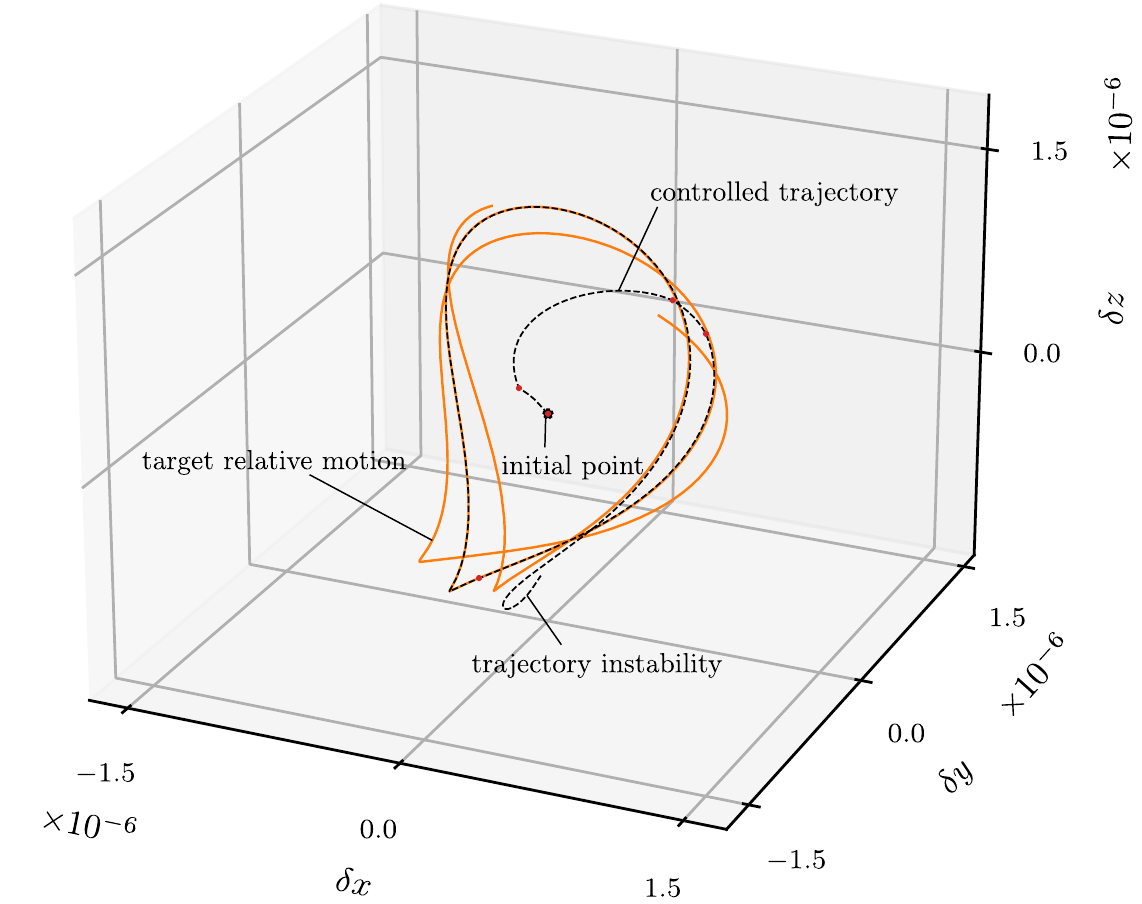}
    \caption{Control to Center Mode, Unstable $L_{2}$ Halo Orbit}%
    \label{fig:HaloControlEx2a}%
\end{figure}
\begin{figure}[h!]
    \centering
    \includegraphics[]{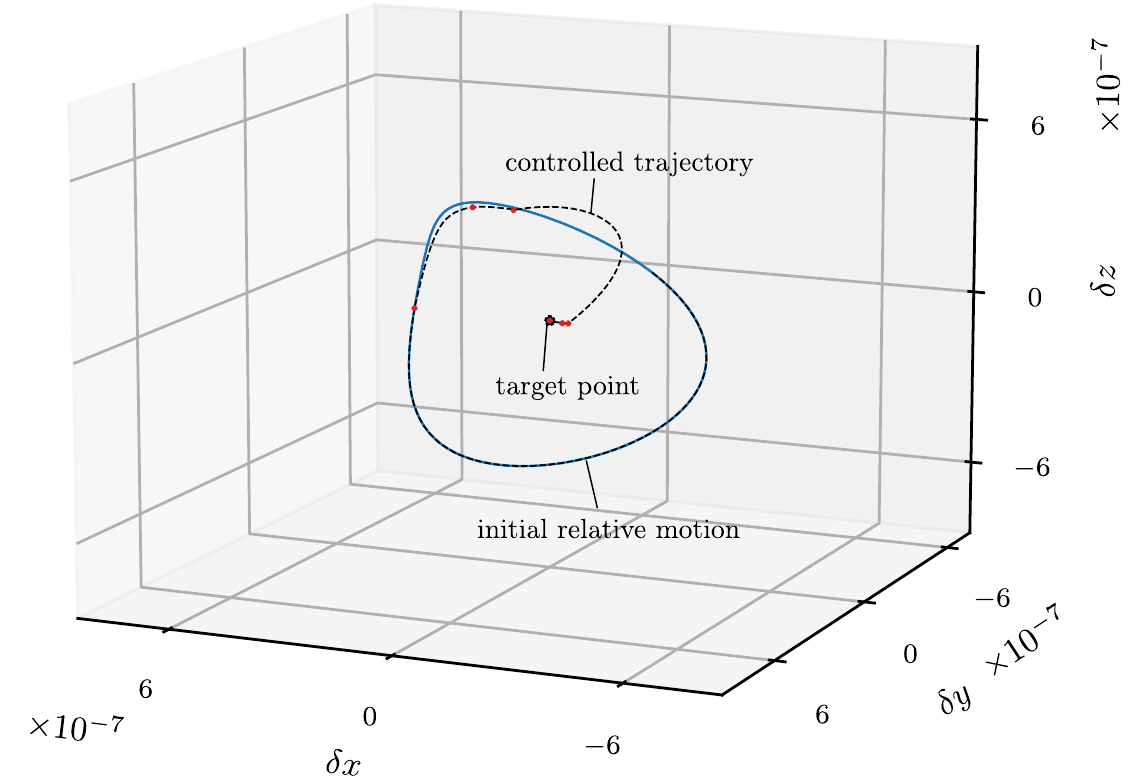}
    \caption{Regulation to Periodic Orbit, Unstable $L_{2}$ Halo Orbit}%
    \label{fig:HaloControlEx2b}%
\end{figure}
\begin{table}[h!]
\centering
\caption{Halo Orbit Control Simulation Parameters, Example 2(a)}
\label{table:L2HaloEx2a}
\begin{tabular}{l|l}
Parameter               & Value                                                                        \\ \hline
Initial relative motion & $\bm{c}_{0} = (0, 0, 0, 0, 0, 0)$ \\
Initial state           & $\delta\bm{x}_{0} = (0, 0, 0, 0, 0, 0)$      \\
Desired relative motion & $\bm{c}^{*} = (0, 0, 0.3, 0, 0, 0)\alpha$, $\alpha = 5.2\times 10^{-6}$                                    \\
Maneuver interval       & $\tau \in [1.90, 5.08]$ \\
Maneuver times & $\tau_{1} = 1.903$, $\tau_{2} = 2.538$, $\tau_{3} = 3.595$, $\tau_{4} = 3.807$, $\tau_{5} = 4.864$ \\
Resulting maneuvers     & $\Delta\bm{v}_{1} = (-0.035, -0.007, 0.052)\alpha$ \\ & $\Delta\bm{v}_{2} = (-0.148, 0.278, 0.052)\alpha$ \\ & $\Delta\bm{v}_{3} = (0.028, 0.002, -0.068)\alpha$ \\ & $\Delta\bm{v}_{4} = (0.028, -0.035, -0.101)\alpha$ \\ & $\Delta\bm{v}_{5} = (-0.001, 0, 0.0005)\alpha$\end{tabular}
\end{table}

\begin{table}[h!]
\centering
\caption{Halo Orbit Control Simulation Parameters, Example 2(b)}
\label{table:L2HaloEx2b}
\begin{tabular}{l|l}
Parameter               & Value                                                                        \\ \hline
Initial relative motion & $\bm{c}_{0} = (0.2, 0, 0, 0, 0, 0)\alpha$, $\alpha = 5.2\times 10^{-6}$ \\
Initial state           & $\delta\bm{x}_{0} = (0., -0.127, 0, -0.117, 0, -0.1)\alpha$      \\
Desired relative motion & $\bm{c}^{*} = (0, 0, 0, 0, 0, 0)$                                    \\
Maneuver interval       & $\tau \in [1.90, 5.08]$ \\
Maneuver times & $\tau_{1} = 1.903$, $\tau_{2} = 2.538$, $\tau_{3} = 2.749$, $\tau_{4} = 4.018$, $\tau_{5} = 4.230$, $\tau_{6} = 4.864$ \\
Resulting maneuvers     & $\Delta\bm{v}_{1} = (-0.008, 0.0005, -0.004)\alpha$ \\ & $\Delta\bm{v}_{2} = (-0.002, -0.003, 0.002)\alpha$ \\ & $\Delta\bm{v}_{3} = (-0.044, 0.084, 0.081)\alpha$ \\ & $\Delta\bm{v}_{4} = (-0.045, -0.086, 0.081)\alpha$ \\ & $\Delta\bm{v}_{5} = (-0.0007, -0.001, 0.0006)\alpha$ \\ & $\Delta\bm{v}_{6} = (-0.008, -0.0004, -0.004)\alpha$ 
\end{tabular}
\end{table}
For this unstable halo orbit, two relative motion control examples are briefly highlighted. First is control from the chief point $\delta\bm{x} = \bm{0}$ to target one of the center modes. The second is an example of regulation from the trivial mode to the chief point. Starting with the example given by Table~\ref{table:L2HaloEx2a} and Figure~\ref{fig:HaloControlEx2a}, the initial point is at $\delta\bm{x} = \bm{0}$ and the target motion is in orange. The controlled trajectory is given by the dashed line, and the control maneuvers are labeled on the plot with red points. The control successfully targets the quasi-periodic mode $\bm{\psi}_{3}$, but small residual error projected into the unstable subspace results in a subsequent need for correction, so the trajectory departs from the target mode. This is a fundamental property of relative motion control in the vicinity of an unstable halo orbit: corrective maneuvers will always be necessary on some timescale, due to the combined effects of nonlinearity and instability.

The second example is given by the data in Table~\ref{table:L2HaloEx2b} and the trajectory in Figure~\ref{fig:HaloControlEx2b}. The initial trivial modal motion is given in blue, the controlled trajectory is given by the dashed line, and the control maneuvers are labeled with red dots. This example demonstrates regulation control in this environment, with the chief at $\delta\bm{x} = \bm{0}$ successfully targeted to a high degree of numerical precision. A similar strategy could be used for orbit regulation, keeping the spacecraft on the unstable periodic orbit. 

\section{Conclusions}
This paper introduces the \textit{method of fundamental modal solutions} for designing, analyzing, and controlling satellite relative motion in the vicinity of general periodic orbits. The close-proximity satellite relative motion is given as a linear sum of fundamental modal solutions, weighed by the modal constants. This facilitates highly simplified geometric exploration of the satellite relative motion problem. Using a variation-of-parameters approach, the modal constants also serve as a rigorous state representation for the perturbed and/or controlled satellite relative state. The representation for perturbed relative motion is tested with the commonly studied $J_{2}$ problem. For the unperturbed problem of satellite relative motion control, the modal constants have null dynamics but are steered in time-varying directions by control. This is a highly beneficial property for control design, especially impulsive control, and this paper also demonstrates the use of modal constants to simplify the control problem using previously established techniques for general integration constant formulations for dynamical systems. 

The procedure applied analytically to the Keplerian relative motion problem can also be applied numerically to other dynamical situations that admit periodic orbits, and this is demonstrated with application to the classical Circular Restricted Three-Body Problem (CR3BP). Highly fuel-efficient impulsive transfer examples are computed using the same procedure as was used for the Keplerian case, in a manner that is also highly numerically efficient. This has many applications. First, it lays the groundwork for computationally efficient on-board close-proximity satellite relative motion transfers in the cislunar environment. The procedure can also be applied in another highly relevant problem, the Augmented Normalized Hill Three-Body Problem (ANH3BP), to explore efficient satellite relative motion control in terminator orbits around asteroids.

The possibilities for future work are extensive. The modal decomposition procedure could be computed in coordinates with superior linearization, yielding a larger region of validity. Additionally, an alternate choice of coordinates for relative motion in the CR3BP could yield greater geometric insight than the standard rotating coordinates used in this work, particularly for the center modes. The current methodology is only applicable to motion in the vicinity of periodic orbits. A perturbative application of the modal decomposition procedure (explored in early work in Reference~\citenum{Burnett2020_Ryugu}) could extend the methods and control techniques discussed in this paper to more realistic setting of motion in the vicinity of a satellite orbit that roughly repeats but is not periodic -- e.g. any orbit for which the chief state $\bm{X}(t) \approx \bm{X}(t+T)$ for some time $T$. This is achieved by finding the nearest $T$-periodic LTV system to a given linearization, and using the modal decomposition from the former to approximately describe the motion in the latter. Additionally, for motion in the vicinity of a general quasi-periodic orbit, the periodic Lyapunov-Floquet reducing transformation is replaced with a quasi-periodic Lyapunov-Perron transformation, which can be challenging to compute. Reference~\citenum{Jorba2001} and references therein are relevant to this. If the reducing transformation for the quasi-periodic case can be computed efficiently, the modal decomposition methodology will be extended quite naturally to general orbits in the full zonal problem and quasi-periodic orbits in the restricted three-body problems. 

\bibliography{references.bib}  

\end{document}